\documentclass[12pt]{amsart}
\usepackage{amsfonts,amssymb}
\usepackage[titletoc]{appendix} 
\usepackage{mathrsfs}
\usepackage{array}
\usepackage[all,ps,cmtip,rotate]{xy} 
\usepackage[usenames,dvipsnames]{color}
\usepackage[hidelinks]{hyperref}
\usepackage{cleveref}
\usepackage{bm,mathtools}
\usepackage{enumerate}
\usepackage{amsmath}
\usepackage{mathrsfs}
\usepackage{geometry}
\usepackage{graphicx}
\usepackage{subfigure}
\usepackage{tikz-cd}
\usepackage{adjustbox}
\usepackage{appendix}
\usepackage{microtype}
\usepackage{geometry}
\newtheorem{define}{Definition}[section]
\newtheorem{thm}[define]{Theorem}
\newtheorem{prop}[define]{Proposition}
\newtheorem{eg}[define]{Example}
\newtheorem{rem}[define]{Remark}
\newtheorem{cor}[define]{Corollary}
\newtheorem{lem}[define]{Lemma}
\newtheorem{conj}[define]{Conjecture}

\def\CH{{\rm CH}}

\def\Z{{\mathbb{Z}}}
\def\e{{\acute{e}}}
\def\nr{\mathrm{nr}}
\def\BM{\mathrm{BM}}
\def\dim{{\mathrm{dim}}}
\def\codim{{\mathrm{codim}}}
\def\P{{\mathbb{P}}}
\def\A{{\mathbb{A}}}
\def\C{{\mathbb{C}}}

\def\Q{{\mathbb{Q}}}
\def\ker{{\mathrm{ker}}}

\def\im{{\mathrm{im}}}
\def\id{{\mathrm{id}}}
\def\cl{{\mathrm{cl}}}

\newcommand{\nospacepunct}[1]{\makebox[0pt][l]{\,#1}}
\date{}
\title[higher Chow groups and refined unramified cohomology]{higher Chow groups with finite coefficients and refined unramified cohomology}
\begin{document}
\author{Kees Kok}
\address{\parbox{0.9\textwidth}{KdV Institute for Mathematics, University of Amsterdam, Netherlands}}
\email{k.kok@uva.nl}
\author{Lin Zhou}
\address{\parbox{0.9\textwidth}{Academy of Mathematics and Systems Science, Chinese Academy of Sciences, Beijing, China;\\ Institute of Algebraic Geometry, Leibniz University Hannover, Hannover, Germany.}}
\email{{linchow@amss.ac.cn; zhou@math.uni-hannover.de}}
\maketitle

\begin{abstract}
In this paper we show that Bloch's higher cycle class map with finite coefficients for quasi-projective equi-dimensional schemes over a field fits naturally in a long exact sequence involving Schreieder's refined unramified cohomology. We also show that the refined unramified cohomology satisfies the localization sequence. Using this we conjecture in the end that refined unramified cohomology is a motivic homology theory and explain how this is related to the aforementioned results.
\end{abstract}
\tableofcontents
\section{Introduction}
In \cite{schreieder2021refined} Schreieder introduced refined unramified cohomology groups, which are a generalization of known theories like unramified cohomology (cf. \cite{colliot1995birational,CTV12}) and Kato homology (cf. \cite{kato1986hasse,KS11}). In particular, Schreieder showed that the refined unramified cohomology groups can detect the failure of the surjectivity and injectivity of the integral cycle class map for cycles of any codimension \cite[Section 7]{schreieder2021refined}, generalizing  \cite{CTV12,voisin2012degree}. If one admits the Hodge/Tate conjecture, this measures the failure of the integral Hodge/Tate conjecture for smooth projective varieties.

A natural question to ask now is, whether the refined unramified cohomology groups can also detect kernels and cokernels of the \emph{higher} cycle class map, defined by Bloch in \cite[Section 4]{bloch11986algebraic}. In other words, we want to compare the differences of motivic homology and Borel--Moore homology in terms of refined unramified cohomology. As the image of the \emph{integral} higher cycle class map over the complex numbers is, for example, always torsion, this might not be the right map to study. It turns out that considering only finite coefficients $M\coloneqq\Z/m\Z$ (here $m$ is invertible in the base field $k$) and the higher cycle class map $\cl^{b,n}\colon\CH^b(X,n;M)\to H_{\e t}^{2b-n}(X,M(b))$ (for smooth $X$) is the natural way to proceed.

It is known that the higher cycle class map $\cl^{b,n}$ with finite coefficients is an isomorphism for $b\leq n$ and injective if $b=n+1$ by the  \emph{Beilinson--Lichtenbaum Conjecture} (see \Cref{cor-GL01}, \cite[Theorem 1.1]{geisser2001bloch} and \cite[Theorem 9.1]{kerz2012cohomological}), moreover, \cite[Theorem 4.2]{suslin2000higher} and \cite[Théorème 1.1]{kelly2014isomorphisme} give that $\cl^{b,n}$ is also an isomorphism for $b\geq\dim(X)$ if the base field $k$ is algebraically closed. The main result of this paper relates the kernel and cokernel of the higher cycle class map in the remaining range of index $b>n+1$ with Schreieder's refined unramified cohomology groups (see \Cref{def-refined unramified cohomology}), denoted by $H^p_{q,\nr}(X,M(n))$, in the form of a long exact sequence:
\begin{thm}[= \Cref{thm-main thm2}]\label{thm-our main intro}
Let $X$ be a $d$-dimensional quasi-projective equi-dimensional scheme defined over a field $k$, and $M\coloneqq \Z/m\Z$, where $m$ is an integer invertible in $k$. For any $b,n\geq 0$, there is a natural long exact sequence
\begin{small}
\begin{align}\label{equ-our main intro}
\to\CH^b(X,n;M)\xrightarrow{\mathrm{cl}}H^{2b-n}(X,M(b))\xrightarrow{\rm restr.} H^{2b-n}_{b-n-1,\nr}(X,M(b))\xrightarrow{\theta} \CH^b(X,n-1;M)\to,
\end{align}  
\end{small}

where $\theta$ is defined in (\ref{equ-construct a connected morphism II}),  $\CH^b(X,n;M)\coloneqq\CH_{d-b}(X,n;M)$ is Bloch's higher Chow group and $H^i(X,M(n))\coloneqq H^{BM}_{2d-i}(X,M(d-n))$ is the $\e$tale Borel--Moore homology with coefficients $M$. The map $\mathrm{cl}\colon\CH^b(X,n;M)\to H^{2b-n}(X,M(b))$ is the higher cycle class map (see \Cref{prop-higher cycle class map for any field}).  
\end{thm}
Here `natural' means that the long exact sequence is compatible with proper pushforward, flat pull-back and pull-back between smooth projective varieties. Note that $X$ in the above statement may be singular. If $X$ is smooth, then the cycle class map $\theta$ appearing in (\ref{equ-our main intro}) agrees with Bloch's cycle class map form \cite{bloch11986algebraic}. The proof of \Cref{thm-our main intro} only uses the Beilinson--Lichtenbaum Conjecture (which also gives the $b\leq n$ case) and we indeed recover \cite{suslin2000higher} again, since \cite[Corollary 5.11]{schreieder2021refined} implies that  $H^{2b-n}_{b-n-1,\nr}(X,M(b))=0$ for $b\geq \dim(X)$ if $\bar{k}=k$.

This result is inspired by a similar result over finite fields involving \emph{Kato homology}, denoted by $KH_*^{\e t}(X,M(n))$, which says that if $k$ is a finite field and $X$ smooth of pure dimension $d$ over $k$, there is a canonical long exact sequence
\[
\cdots\to \CH^d(X,n;M)\xrightarrow{\cl} H_{\e t}^{2d-n}(X,M(d))\to KH_{n+1}^{\e t}(X,M)\to \CH^d(X,n-1;M)\to \cdots,
\]
proven in \cite[Lemma 6.2]{jannsen2009kato}. 
Motivated by the above exact sequence for smooth schemes, we believe that the exact sequence (\ref{equ-our main intro}) should also hold for arbitrary fields and schemes.

A direct consequence of \Cref{thm-our main intro} is that the kernel of the natural restriction map $H^{b+1}(X,M(b))\to H^{b+1}_{\nr}(X,M(b))$ is isomorphic to $\CH^b(X,b-1;M)$ for smooth quasi-projective $X$ (i.e., \Cref{cor-kernel of map}). Moreover, we obtain an exact sequence related to Bloch's map on torsion cycles.
\begin{cor}[= {\Cref{cor-an exact sequence about torsion cycle}}]
Let $X$ be a smooth quasi-projective variety defined over an algebraically closed field $k$, and let $\ell$ be a prime which is not equal to $\mathrm{char}(k)$. There is an exact sequence
\begin{equation}
H^{2n-2}_{n-3,\nr}(X,\Q_{\ell}/\Z_{\ell}(n))\to \CH^n(X)[\ell^{\infty}]\xrightarrow{\lambda^n} H^{2n-1}(X,\Q_{\ell}/\Z_{\ell}(n))\to H^{2n-1}_{n-2,\nr}(X,\Q_{\ell}/\Z_{\ell}(n)),
\end{equation}
where $\lambda^n$ is the morphism defined by Bloch in \cite{bloch1979torsion} and the first arrow equals the composition $ H^{2n-2}_{n-3,\nr}(X,\Q_{\ell}/\Z_{\ell}(n))\to \CH^n(X,1;\Q_{\ell}/\Z_{\ell})\to \CH^n(X)[\ell^{\infty}]$. If $k=\C$, we can replace $\mathbb{Q}_{\ell}/\mathbb{Z}_{\ell}$ (resp. $\CH^n(X)[\ell^{\infty}]$) by $\mathbb{Q}/\Z$ (resp. $\CH^n(X)_{\rm tor}$).
\end{cor}

Another result of this paper comes from the observation that `two out of three' of the groups involved in the exact sequence of \Cref{thm-our main intro} satisfy certain properties of a cohomology theory. It turns out that the refined unramified cohomology is indeed a homology theory with properties similar to Bloch's higher Chow groups and Borel--Moore homology: 
\begin{thm}[\Cref{properties of homology theory}]\label{thm-refined unramified homology is homology}
Let $X$ be a separated scheme of finite type over a field $k$. The refined unramified cohomology $H^p_{q,\nr}(X,A(n))$ is a homology theory satisfying the following conditions
\begin{itemize}
\item (Localization Sequence) Let $i\colon Z\to X$ of be a closed immersion of codimension $c\coloneqq\dim(X)-\dim(Z)$ with open complement $j\colon U\hookrightarrow X$ such that $\dim(U)=\dim(X)$. There is a natural long exact sequence
\begin{equation*}
    \cdots\to H^{p-2c}_{q-c,\nr}(Z,A(n-c))\xrightarrow{i_*}H^{p}_{q,\nr}(X,A(n))\xrightarrow{j^*}H^{p}_{q,\nr}(U,A(n))\xrightarrow{\delta_{q}^{p}}H^{p
    +1-2c}_{q+1-c,\nr}(Z,A(n-c))\xrightarrow{i_\ast}\cdots.
\end{equation*}
\item ($\mathbb{A}^1$-Homotopy Invariance) Let $f\colon Y\to X$ be an affine bundle. Suppose the twisted Borel-Moore cohomology $H^i(-,A(n))$ is $\mathbb{A}^1$-homotopy invariant (see \Cref{def-A1 homotopy invariance}). The pullback $f^*\colon H^p_{q,\nr}(X,A(n))\to H^p_{q,\nr}(Y,A(n))$ induced by the projection is a natural isomorphism for all $p$ and $q$.

\end{itemize}
\begin{itemize}
\item (Bockstein Sequence) If the twisted Borel--Moore cohomology $H^i(-,A(n))$ satisfies the Bloch--Kato $\ell$-type conjecture (see \Cref{def-BK type}), then there is a natural morphism $H^p_{q,\nr}(X,\mu_{\ell^r}^{\otimes (p-q-1)})\xrightarrow{\Delta^p_q}H^{p+1}_{q+1,\nr}(X,\Z_{\ell}(p-q-1))$ for any integers $r$, $p$ and $q$, such that 
\begin{small}
\begin{equation*}
   \to H^p_{q,\nr}(X,\Z_\ell(p-q-1))\xrightarrow{\times \ell^r}H^p_{q,\nr}(X,\Z_\ell(p-q-1))\to  H^p_{q,\nr}(X,\mu_{\ell^r}^{\otimes (p-q-1)})\xrightarrow{\Delta^p_q} H^{p+1}_{q+1,\nr}(X,\Z_\ell(p-q-1))\xrightarrow{\times \ell^r}     
\end{equation*}
\end{small}
\noindent is natural exact. Furthermore, if $H^i(-,A(n))$ is the pro-\'etale cohomology (cf.\cite{BS15}) and  $k$ contains all $\ell$-power roots of unity, there exists a long exact sequence:
\[
\cdots\to H^p_{q,\nr}(X,\Z_{\ell}(n))\xrightarrow{\times \ell^r} H^p_{q,\nr}(X,\Z_{\ell}(n))\to H^p_{q,\nr}(X,\mu_{\ell^r}^{\otimes n})\xrightarrow{\Delta^p_q}H^{p+1}_{q+1,\nr}(X,\Z_{\ell}(n))\xrightarrow{\times \ell^r}\cdots.  
\]
If $k=\C$, then we can replace $\Z_\ell$ and $\mu_{\ell^r}$ by $\Z$ and $\Z/n\Z$ for any $n$, respectively.
\end{itemize}
\end{thm}
After introducing the main definitions in \Cref{setion-homology theory}, we shall prove these properties in \Cref{properties of homology theory}. As an application of the localization sequence and the homotopy invariance, we also give a computation of the refined unramified cohomology of the trivial projective bundle in \Cref{computation trivial proj bundle} (cf.~ \cite[Corollary 1.8]{schreieder2022moving} for the similar result involving smooth projective equi-dimensional $k$-schemes). For a computation of general projective bundles and a blow-up formula, we refer to \cite{KZ23functorial}. 

Moreover, the localization sequence also induces a `coniveau filtration' spectral sequence as in \cite{BO74}, converging to the refined unramified cohomology (see (\ref{equ-coniveau for refined})). In the end, we conjecture that the refined unramified cohomology has the nature of a motivic homology theory defined over a certain motivic category of any field. Specifically, we conjecture it is a sheaf cohomology for smooth $X$. During the reviewing process of the present paper, this conjecture has been extended and is proven by Alexandrou and Schreieder (see \cite{alexandrou2024truncated}).
\begin{conj}[= \Cref{conj-truncation for smooth case}]\label{conj}
Let $A(n)$ be a (twisted) torsion abelian group whose elements have order coprime to the exponential characteristic of $k$, and let $X$ be a smooth variety defined over a field $k$. Consider the identity continuous morphism $\pi\colon X_{\rm \e t}\to X_{\rm Zar}$, and let $\tau^{\geq l}R\pi_*A(n)$ be the $l$-th canonical truncation of the derived complex $R\pi_*A(n)$, where $A(n)$ is the locally constant sheaf with stalk $A(n)$. There is a canonical isomorphism $\mathbb{H}^k(X,\tau^{\geq l}R\pi_*A(n))\cong H^k_{k-l,\nr}(X,A(n))$ for any $k,l\geq 0$. Moreover, for the case of $F=\mathbb{C}$, $A$ can be any abelian group.
\end{conj}

\subsection{Notations and Conventions}
\begin{itemize}
   \item  We let $k$ be a field and $m$ an integer invertible in $k$. We write $M\coloneqq \Z/m$ without emphasis. 
   \item  All schemes are separated and Noetherian, generally denoted by $X$. 
   \item An algebraic scheme is a separated scheme of finite type over a field. 

   \item For a scheme $X$, let $X^{(p)}=\{x\in X\mid {\codim}(x)=p\}$ be the set of $p$-codimensional points in $X$, where $\codim (x):=\dim(X)-\dim (x)$. Let $F_pX=\{x\in X|\codim(x)\leq p\}$. Similarly, let $X_{(p)}=\{x\in X\mid {\dim}(x)=p\}$ be the set of $p$-dimensional points in $X$, and let $F^pX=\{x\in X|\dim(x)\geq p\}$.
   \item We often denote the higher cycle class map $\cl^{i,n}_M$ (or $\cl_{i,n}^M$) by `$\cl$' for readability reasons.
   \item  We always assume that smooth schemes are equi-dimensional.
\end{itemize}

\subsection{Acknowledgements}
The authors are very grateful to Simon Pepin Lehalleur for his help in finding and explaining the references used in \Cref{construction higher cycle class map}.  The authors also would like to thank the reviewer for his/her valuable comments which greatly enhance the writing of this paper. The first author would also like to thank Mingmin Shen and Rob de Jeu for the helpful conversations. The second author would like to thank Zhiyu Tian and Enlin Yang for their useful suggestions. She also thanks Stefan Schreieder for his discussions during the revision of the paper. Simultaneously, she has received funding from the European Research Council under the European Union's Horizon 2020 research and innovation program under grant agreement
No 948066 (ERC-StG RationAlgic).

\section{Homology Theory}\label{setion-homology theory}
In this section, we recall some content about homology theories and refined unramified cohomology for convenience of the reader. Let $\mathcal{C}$ be a category of algebraic $k$-schemes containing all quasi-projective  $k$-schemes such that every closed (or open) immersion is a morphism in $\mathcal{C}$. Thus every locally closed algebraic $k$-subscheme of an object $X$ in $\mathcal{C}$ is also an object in $\mathcal{C}$. 
\begin{define}[{\cite[Definition 2.1]{JS03}}]\label{def-homology theory}
Let $\mathcal{C}_*$ be the category whose objects are the same as $\mathcal{C}$, and the sets of morphisms consist only of proper morphisms. A homology theory on $\mathcal{C}$ is a sequence of covariant functors ($i\in \mathbb{Z}$)
\begin{equation}
H_i\colon \mathcal{C}_*\to \{\text{category of abelian groups}\}
\end{equation}
satisfying the following conditions:
\begin{enumerate}	
\item[(a)]  For any open immersion $j\colon V\hookrightarrow X$ in $\mathcal{C}$, there is a functorial map $j^*\colon H_i(X)\to H_i(V)$.
\item[(b)] If $i\colon Y\hookrightarrow X$ is a closed immersion of $X$ with open complement $j\colon V\hookrightarrow X$, there is a long exact sequence (called localization sequence)
\begin{equation*}
\cdots\xrightarrow{\partial}H_i(Y)\xrightarrow{i_*}H_i(X)\xrightarrow{j^*}H_i(V)\xrightarrow{\partial}H_{i-1}(Y)\to\cdots.
\end{equation*}
The maps $\partial$ are called the connecting morphisms. This sequence is functorial with respect to proper maps or open immersions in a natural way. Moreover, a morphism between homology theories $H$ and $H'$ is a morphism $\phi\colon H\to H'$ of functors on $\mathcal{C}$ compatible with the localization sequence in the natural way.
\end{enumerate}
\end{define}
\begin{eg}
The \'etale Borel--Moore homology theory with finite coefficients (cf. \cite[Section 2]{BO74}) and Bloch's higher Chow groups (cf. \cite[Section 3]{bloch1986algebraic}) are examples of homology theories.
\end{eg}
For a fixed homology theory $H_*$ on $\mathcal{C}$, and for every object $X$ in $\mathcal{C}$, Bloch and Ogus constructed a spectral sequence associated with the \emph{niveau filtration}; see \cite{BO74}.
\begin{prop}[{\cite[Proposition~3.7]{BO74}}]\label{prop-niveau filtration}
If $H_*$ is a  homology theory on $\mathcal{C}$, for any object $X$ in $\mathcal{C}$, there is a spectral sequence 
\begin{equation}\label{equ-niveau filtration}
E^1_{p,q}(X)=\bigoplus_{x\in X_{(p)}}H_{p+q}(x)\Rightarrow N_{\bullet}H_{p+q}(X),
\end{equation}
where  
\begin{equation*}
H_a(x)\coloneqq \varinjlim_{V\subset\overline{\{x\}}}H_a(V),
\end{equation*}
the direct limit is taken over all nonempty open subschemes of $\overline{\{x\}}$, and 
\begin{equation*}
N_pH_{p+q}(X)
= \bigcup_{Z\subset X}\mathrm{Im}\{H_{p+q}(Z)\to H_{p+q}(X)\},
\end{equation*}
where $Z$ is taken over all closed subsets of $X$ with $\dim(Z)\leq p$.
This spectral sequence is covariant with respect to proper morphisms in $\mathcal{C}$ and contravariant with respect to open immersions. For every morphism $\phi\colon H\to H'$ between homology theories, $\phi$ induces a morphism between the associated niveau spectral sequences.    
\end{prop}
For every fixed $q$ and a fixed homology theory $H_*$ over $\mathcal{C}$, the group $E^1_{p,q}(X)$ fits in a short exact sequence
\begin{equation*}
   0\to E^1_{p,q}(Y)\xrightarrow{i_*} E^1_{p,q}(X)\xrightarrow{j^*} E^1_{p,q}(U)\to 0 
\end{equation*}
by \cite[Proposition~2.9]{JS03} and \cite[Proposition~7.2, Proposition~7.3]{KS11}, where $i:Y\to X$ is a closed subscheme, $j:U\to X$ is the open complement of $Y$ and $i_*$ (resp. $j^*$) is induced by proper pushforward (resp. open immersion). This short exact sequence in turn induces the following long exact sequence
\begin{equation*}
   \cdots \to E^2_{p,q}(Y)\to E^2_{p,q}(X)\to E^2_{p,q}(U)\to E^2_{p-1,q}(Y)\to \cdots.
\end{equation*}
In others words, for every fixed $q$, the family of functors $\{E^2_{p,q}\}_{p\in \mathbb{Z}}$ defines a homology theory on $\mathcal{C}$; see \cite[Corollary 2.10]{JS03}.

Similarly, one can also define a twisted homology theory $H_i(-,n)$ and construct the twisted spectral sequence associated with the niveau filtration as \cite[Definition~1.2, Proposition~3.7]{BO74}. Moreover, similarly for a twisted homology theory $H_i(-,n)$ over $\mathcal{C}$, the family of functors $\{E^2_{p,q}\}_{p\in \mathbb{Z}}$ also defines a twisted homology theory on $\mathcal{C}$ for every fixed $q$.
\subsection{Borel--Moore Homology and Refined Unramified Cohomology}\label{section-2.1}
Let $X$ be an algebraic scheme defined over a field $k$ of dimension $d$ and let $\pi_X\colon X\to{\rm Spec}(k)$ be the structure map of $X$. Let $\ell$ be a prime invertible in $k$ and let $m$ be an integer in $k^*$. For every $n\in \Z$, we consider the \emph{pro-étale cohomology} (see \cite{BS15}):
\begin{itemize}
    \item[(a)]\label{eg-a} $H^i(X,\mu_m^{\otimes n})\coloneqq R^{i-2d}\Gamma(X_{\text{proét}},\pi_X^!\mu_m^{\otimes n-d})$;
    \item[(b)] $H^i(X,\Z_{\ell}(n))\coloneqq R^{i-2d}\Gamma(X_{\text{proét}},\pi_X^!\Hat{\Z}_{\ell}(n-d))$;
    \item[(c)] $H^i(X,\Q_{\ell}/\Z_{\ell}(n))\coloneqq \varinjlim\limits_r H^i(X,\mu_{\ell^r}^{\otimes n})$, \quad $H^i(X,\Q_{\ell}(n))\coloneqq H^i(X,\Z_{\ell}(n))\otimes_{\Z_{\ell}}\Q_{\ell}$.
\end{itemize}
\begin{rem}
\rm{By the comparison theorem \cite[Proposition 5.2.6]{BS15}, for all $i,n$ and $m$, there is a natural isomorphism 
\begin{equation}\label{equ-comparison theorem}
  H_{2d-i}^{BM}(X,\mu_m^{\otimes d-n})\coloneqq  H^{i-2d}(X_{\text{ét}},\pi_X^!\mu_m^{\otimes n-d})\cong  R^{i-2d}\Gamma(X_{\text{proét}},\pi_X^!\mu_m^{\otimes n-d})\eqqcolon H^i(X,\mu_m^{\otimes n}),
\end{equation}
  where the first group is exactly the \emph{$\e$tale Borel--Moore homology} with finite coefficients. Therefore, if $k$ is either algebraically closed or finite, we have an isomorphism $$H^i(X,\Z_{\ell}(n))\cong H_{2d-i}^{BM}(X_{\text{ét}},\Z_{\ell}(d-n))\coloneqq \varprojlim\limits_r H_{2d-i}^{BM}(X_{\text{ét}},\mu_{\ell^r}^{\otimes d-n})$$ for every $n$ and $i$. Moreover for $k=\C$, we also consider the \emph{Borel-Moore singular homology} for the abelian group A:
 \begin{itemize}
     \item[(d)] $ H^i(X,A(n) )\coloneqq H^{\BM}_{2d_X-i}(X_{\rm an},A).$
\end{itemize}
Note that if $A\in \{\mu_{m},\mathbb Z_\ell,\Q_{\ell}/\Z_{\ell},\Q_{\ell}\}$, (d) coincides with the pro-étale cohomology groups (étale cohomology groups in this case) defined in (a)-(c) by another comparison theorem (cf. \cite[Theorem 21.1]{milne2012lectures}). Thus the above definitions are compatible.
}
\end{rem}
The above mentioned are examples of  \emph{twisted Borel--Moore cohomologies} (see \cite[Section~4 and Section~5]{schreieder2021refined}) with flat pull-backs.

\textbf{Notation:}
For any algebraic $k$-scheme $X$ of dimension $d$, suppose $$A(n)\in \{\mu_{m}^{\otimes n},\mathbb Z_\ell(n),\Q_{\ell}/\Z_{\ell}(n),\Q_{\ell}(n)\}.$$ We write $H^i(X,A(n))$ to represent the groups (a)-(c) without emphasis. For $k=\C$, we also write $H^i(X,A(n))$ to represent the groups (d) for the abelian group $A$.  Moreover, we sometimes use the homological notation  $H_i(X,A(n))\coloneqq H^{2d-i}(X,A(d-n))$ for convenience.

The group \( H_*(X, A(n)) \) is a (twisted) homology theory satisfying \Cref{def-homology theory}. We also have the \emph{Gysin exact sequence} (cf. \Cref{prop-formal Gysin exact sequence}), which can be considered as the localization sequence, leading to the existence of a niveau spectral sequence associated with \( H^i(X, A(n)) \). In particular, for a smooth \( X \) of pure dimension, \( H^i(X, A(n)) \) is canonically isomorphic to ordinary cohomology by \emph{Poincar\'e duality}. Furthermore, the ordinary cohomology of the aforementioned cases (a)-(d) are examples of \emph{a twisted cohomology theory with an action by cycles} (see \cite[Section~3]{schreieder2022moving}).

\begin{prop}[{\cite[Definition 4.2, Proposition 6.6 and 6.9]{schreieder2021refined}}]\label{prop-formal Gysin exact sequence}
Let $X$ be a Noetherian scheme and let $i\colon Z\hookrightarrow X$ be a closed immersion of codimension $c\coloneqq\dim(X)-\dim (Z)$ and with complement $j\colon U\to X$ with $\dim(X)=\dim(U)$. Suppose $H^*(-,A(n))$ is a twisted Borel–Moore cohomology theory. Then there is a Gysin exact sequence
\begin{equation}\label{equ-formal gysin sequence}
    \cdot\cdot\to H^i(X,A(n))\to H^i(U,A(n))\xrightarrow{\partial}H^{i+1-2c}(Z,A(n-c))\xrightarrow{\iota_*}H^{i+1}(X,A(n))\to\cdot \cdot , 
\end{equation}
where $\partial$ is called the residue map.
\end{prop}
\begin{define}[{\cite[Definition 5.1]{schreieder2021refined}}]\label{def-refined unramified cohomology}
For any Noetherian scheme $X$, let $$H^i(F_jX,A(n))\coloneqq \varinjlim_{F_jX\subset U\subset X}H^i(U,A(n)).$$ The $j$-th refined unramified cohomology of $X$ with coefficient $A(n)$ is defined as 
\[H^i_{j,\nr}(X,A(n))\coloneqq \im (H^i(F_{j+1}X,A(n))\to H^i(F_jX,A(n))).\]
\end{define}
Similarly, we can define the refined unramified homology.
\begin{define}\label{def-refined umramified homology}
 For any Noetherian scheme $X$, let $H_i(F^jX,A(n))\coloneqq\varinjlim_{F^jX\subset U\subset X}H_i(U,A(n))$. The $j$-th refined unramified homology of $X$ with coefficient $A(n)$ is defined as 
\[H_i^{j,\nr}(X,A(n))\coloneqq\im (H_i(F^{j-1}X,A(n))\to H_i(F^jX,A(n))).\]   
\end{define}
By this definition, we have $H^{d-j,\nr}_{2d-i}(X,A(d-n))=H^i_{j,\nr}(X,A(n))$ where $d\coloneqq\dim(X)$. Taking direct limits of the sequence in \Cref{prop-formal Gysin exact sequence} twice, we get a \emph{limit version of the Gysin sequence}.
\begin{prop}[{\cite[Lemma 5.7]{schreieder2021refined}}]\label{prop-Gysin sequence}
Let $X$ be a Noetherian scheme; then for any $m,j\geq 0$ there is a canonical long exact sequence
\[
\cdots\to H^i(F_{j+m}X,A(n))\to H^i(F_{j-1}X,A(n))\overset{\partial}{\to}\varinjlim_{Z\subset X} H^{i+1-2j}(F_mZ,A(n-j))\overset{\iota_\ast}{\to}\cdots,
\]
where $Z$ runs through all closed subsets of $X$ with $\codim(Z)\coloneqq \dim(X)-\dim(Z)=j$.
\end{prop}
If we take $m=0$ in \Cref{prop-Gysin sequence} then we obtain an exact couple $$D_1\xrightarrow{(-1,1)}D_1\xrightarrow{(0,0)}E_1\xrightarrow{(1,0)}D_1$$ with $D_1^{p,q}\coloneqq H^{p+q-1}(F_{p-1}X,A(n))$ and $E_1^{p,q}(X,A(n))\coloneqq\bigoplus\limits_{x\in X^{(p)}}H^{q-p}(x,A(n-p))$; see \cite[\textsection 1.3]{schreieder2021refined}. The next page of this spectral sequence then gives the following.

\begin{prop}[{\cite[Proposition 7.35]{schreieder2021refined}}]\label{prop-exact sequence of D2}
Let $X$ be a Noetherian scheme. There is a canonical long exact sequence
\[
\cdots\to H_{j-1,\nr}^{i+2j-1}(X,A(n))\to H_{j-2,\nr}^{i+2j-1}(X,A(n))\to E_2^{j,i+j}(X,A(n))\to H_{j,\nr}^{i+2j}(X,A(n))\to\cdots,
\]
where 
\begin{equation*}
    E_2^{p,q}(X,A(n))\coloneqq \frac{\ker(E_1^{p,q}(X,A(n))\xrightarrow{\partial\circ \iota_*}E_1^{p+1,q}(X,A(n)))}{\im(E_1^{p-1,q}(X,A(n))\xrightarrow{\partial\circ \iota_*} E_1^{p,q}(X,A(n)))}.
\end{equation*}
\end{prop}
\begin{cor}[{\cite[Corollary 5.10]{schreieder2021refined}}]\label{cor-degeneration of refined unramified cohomology}
 For any Noetherian scheme $X$, for all $j\geq \lceil i/2\rceil$ we have $H^i(F_jX,A(n))\cong H^i(X,A(n))$.
\end{cor}
\begin{rem}
\rm{ We sometimes use the homological notations $H_*(X,A(n))$ and $H_*^{\dagger,\nr}(X,A(n))$. Then \Cref{prop-Gysin sequence} and \Cref{prop-exact sequence of D2} say that there exist canonical long exact sequences 
\begin{equation}\label{equ-exact sequences for homology}
 \cdot \cdot\to H_i(F^{j-m}X,A(n))\xrightarrow{\rm restr.} H_i(F^{j+1}X,A(n))\xrightarrow{\partial}\varinjlim_{\substack{Z\subset X\\
\dim(Z)=j}} H_{i-1}(F^{j-m}Z,A(n))\overset{\iota_\ast}{\to}\cdot\cdot;  
\end{equation}
\begin{equation}\label{equ-exact sequences for refined unramified homology}
\cdot \cdot\to H^{p+1,\nr}_{p+q+1}(X,A(n))\to H^{p+2,\nr}_{p+q+1}(X,A(n))\to E^2_{p,q}(X,A(n))\to H^{p,\nr}_{p+q}(X,A(n))\to\cdot \cdot,
\end{equation}
where the direct limit runs through all closed subsets $Z\subset X$ with $\dim(Z)=j$.
    }
\end{rem}
The aforementioned properties and definitions are valid for the general twisted Borel–Moore cohomology theory  $H^*(-,A(n))$ on constructible categories (see \cite[Definition 4.1]{schreieder2021refined}).
\begin{thm}[{\cite[Theorem 6.1]{voevodsky2011motivic}}]\label{thm-Bloch Kato conjecture}
Let $k$ be a field which contains $1/\ell$. Then the norm residue homomorphisms
$$K_n^M(k)/{\ell}\to H_{\e t}^n(k,\mu_{\ell}^{\otimes n})$$
are isomorphisms for all $n$.
\end{thm}
Applying the \emph{Bloch--Kato conjecture} (i.e., \Cref{thm-Bloch Kato conjecture}), there is an important corollary by Bloch; see \cite[end of Lecture 5]{bloch2010lectures}.

\begin{cor}\label{cor-torsion free}
Let $X$ be an algebraic $k$-scheme and $x\in X$. The group $H^n(x,\Z_\ell(n-1))$ is torsion free for all $n\geq 0$. If $k=\C$, then group $H^n(x,\Z)$ is torsion free for all $n\geq 0$.
\end{cor}
\subsection{Functoriality of Refined Unramified Cohomology}\label{section-functoriality}
In this section we recall that refined unramified cohomology is motivic and satisfies various natural functoriality properties, thanks to the
works of Schreieder; see \cite{schreieder2020infinite,schreieder2022moving}. 

\begin{prop}\label{prop-functorial pushforward}
Let $f\colon X\to Y$ be a proper morphism with  relative codimension $c\coloneqq \dim(Y)-\dim(X)$. There exists a functorial pushforward map $f_\ast\colon H^i(F_jX,A(n))\to H^{i+2c}(F_{j+c}Y,A(n+c))$. 
If $f\colon X\to Y$ is flat and the twisted Borel--Moore cohomology theory $H^i(-,A(n))$ has flat pull-backs, then there exists a functorial pullback map $f^\ast\colon H^i(F_jY,A(n))\to H^i(F_jX,A(n))$.
\end{prop}
\begin{proof}
This is done in \cite[Lemma 2.5]{schreieder2020infinite} for the étale and continuous cohomology theories. We note here that the same proof works for any twisted Borel--Moore cohomology theory having a flat pull-back and proper push-forward. In particular, the result holds for the examples (a)-(d) from \Cref{section-2.1}. 
\end{proof}
\begin{cor}\label{cor-pushforward for proper}
Let $f\colon X\to Y$ be a proper morphism with  relative codimension $c\coloneqq \dim(Y)-\dim(X)$. There exists a functorial pushforward map $f_\ast\colon H^i_{j,\nr}(X,A(n))\to H^{i+2c}_{j+c,\nr}(Y,A(n+c))$. If $f$ is flat and the twisted Borel--Moore cohomology theory $H^i(-,A(n))$ has flat pull-backs, then there is a functorial pullback $f^\ast\colon H^i_{j,\nr}(Y,A(n))\to H^i_{j,\nr}(X,A(n))$.
\end{cor}
For a general twisted Borel--Moore cohomology theory $H^*(-,A(n))$, the proper pushforward between refined unramified cohomology is also well-defined. However, flat pullbacks are a bit more subtle, as \( H^*(-, A(n)) \) itself may not have flat pullbacks.

Moreover, the refined unramified cohomology associated to the examples (a)-(d) in \Cref{section-2.1} is motivic, which means that it has action of correspondences for smooth projective varieties and has pullbacks between smooth schemes.
\begin{prop}[{\cite[Corollary 6.9]{schreieder2022moving}}]\label{prop-functorial pullback}
If $f\colon X\to Y$ is a morphism between smooth quasi-projective equi-dimensional schemes with smooth compactification (e.g. smooth quasi-projective varieties over ${\rm char}(k)=0$), then there is a functorial pullback morphism $f^\ast\colon H^i_{j,\nr}(Y,A(n))\to H^i_{j,\nr}(X,A(n))$.
\end{prop}
\begin{prop}[{\cite[Corollary 1.7]{schreieder2022moving}}]\label{prop-correspondences actions}
 Let $X$ and $Y$ be smooth projective equi-dimensional schemes over a field $k$ with $d_X=\dim(X)$. For any $c,i,j\geq 0$, there is a natural bi-additive pairing
 $$\CH^c(X\times Y)\times H^i_{j,\nr}(X,A(n))\to H^{i+2c-2d_X}_{j+c-d_X,\nr}(Y,A(n+c-d_X)),$$
 which is functorial with respect to the composition of correspondences.
\end{prop} 

\section{Refined Unramified Cohomology is a Homology Theory}\label{properties of homology theory}
In this section, we show that the refined unramified cohomology has certain properties expected from a homology theory, most notably, it has a localization sequence. As an application, we compute the refined unramified cohomology of trivial projective bundles. 
\subsection{Localization Sequence}
Here we show that there is a localization exact sequence for refined unramified cohomology theories. We use the homological notation from \Cref{def-refined umramified homology} for convenience.

\begin{lem}\label{lem-construct delta}
Let $X$ be a Noetherian scheme and let $Z\subset X$ be a closed subscheme with open complement $U$. There is a natural morphism $\delta_i^j\colon H_i^{j,\nr}(U,A(n))\to H_{i-1}^{j-1,\nr}(Z,A(n))$  for any $i$ and $j$, induced by the residue map for $j\leq \lfloor i/2\rfloor$. 
\end{lem}
\begin{proof}
First we show that there is a morphism $\delta\colon H_i(F^{j-1}U,A(n))\to H_{i-1}(F^{j-2}Z,A(n))$ for any $i,j$, such that the diagram 
\begin{equation}\label{diagram-delta}
\begin{tikzcd}
 H_i(F^{j-1}U,A(n))\ar[r,"\delta"]\ar[d,"{\rm restr.}"]& H_{i-1}(F^{j-2}Z,A(n))\ar[d,"{\rm restr.}"]\\
H_i(F^{j}U,A(n))\ar[r,"\delta"]& H_{i-1}(F^{j-1}Z,A(n))
\end{tikzcd}
\end{equation}
commutes. This in turn defines the desired morphism $\delta_i^j\colon H_i^{j,\nr}(U,A(n))\to H_{i-1}^{j-1,\nr}(Z,A(n))$.

For any $\alpha\in  H_i(F^{j-1}U,A(n))$, we can find an open subset $V\subseteq U\subseteq X$ such that $\dim(U\setminus V)=j-2$ and $\alpha$ can be lifted to $ H_i(V,A(n))$. Let $W$ be the closure of $U\setminus V$ in $X$. Then $\dim(W)=\dim (U\setminus V)=j-2$ and $\dim(W\cap Z)\leq j-3$. Consider the closed immersion $Z\setminus(W\cap Z)\subseteq X\setminus W$ with complement $V$. The residue map  $\partial\colon  H_i(V,A(n))\to  H_{i-1}(Z\setminus W,A(n))$ induced by the Gysin exact sequence, induces the morphism $\delta\colon  H_i(F_{j-1}U,A(n))\to  H_{i-1}(F_{j-2}Z,A(n))$. Both, the fact that $\delta$ is well-defined and that it makes (\ref{diagram-delta}) commute, follow from the fact that the residue map is compatible with restrictions to opens. More precisely, if $V'\subseteq V$ are opens with corresponding $W\subseteq W'$ constructed as above, then we have a commutative diagram
\[
\begin{tikzcd}
H_i(V,A(n))\ar[r,"\partial"]\ar[d,"{\rm restr.}"]& H_{i-1}(Z\setminus W,A(n))\ar[d,"{\rm restr.}"]\\
    H_i(V',A(n))\ar[r,"\partial"]& H_{i-1}(Z\setminus W',A(n)).
\end{tikzcd}
\]

Finally, the naturality of $\delta_i^{j}$ follows from the naturality properties of the residue map on ordinary Borel--Moore cohomology. For the cases $j\leq \lfloor i/2\rfloor$, we apply \Cref{cor-degeneration of refined unramified cohomology}.
\end{proof}
\begin{prop}\label{prop-localization sequence}
Let $X$ be a Noetherian scheme and let $i\colon Z\to X$ be a closed immersion with open complement $j\colon U\hookrightarrow X$. There is a canonical exact sequence
\begin{equation}\label{equ-localization sequence for refined homology}
\to H_{p}^{q,\nr}(Z,A(n))\xrightarrow{i_*}H_{p}^{q,\nr}(X,A(n))\xrightarrow{j^*}H_{p}^{q,\nr}(U,A(n))\xrightarrow{\delta^{q}_{p}}H_{p-1}^{q-1,\nr}(Z,A(n))\to,
\end{equation}
where $i_*$ and $j^*$ is induced by the functoriality of $H_*$. In particular, $H_p^{q,\nr}(X,A(n))$ is a homology theory.
\end{prop}
\begin{proof}
We first show the exactness of the following sequence
\begin{small}
\begin{equation}\label{equ-exact sequence of F^q}
    \to H_p(F^{q}Z,A(n))\xrightarrow{i_*} H_p(F^{q}X,A(n))\xrightarrow{j^*} H_p(F^{q}U,A(n))\xrightarrow{\rm restr.\circ \delta} H_{p-1}(F^qZ,A(n))\to.
\end{equation}
\end{small}

Let $W\subset X$ be a closed subset with $\dim(W)=q-1$. Then we have the Gysin exact sequence
\begin{equation*}
    \to H_p(Z\setminus W,A(n))\to H_p(X\setminus W,A(n))\to H_p(U\setminus W,A(n))\xrightarrow{\partial}H_{p-1}(Z\setminus W,A(n))\to.
\end{equation*}
Taking the direct limits of $W$ in the exact sequence above, we obtain the desired exact sequence (\ref{equ-exact sequence of F^q}).

 To prove (\ref{equ-localization sequence for refined homology}), we first show the sequence (\ref{equ-localization sequence for refined homology}) is exact at $H_p^{q,\nr}(X,A(n))$. By (\ref{equ-exact sequence of F^q}), it remains to show $\ker(j^*)\subset \im(i_*)$, and if $$\alpha\in\ker(H_{p}^{q,\nr}(X,A(n))\xrightarrow{j^*}H_{p}^{q,\nr}(U,A(n)))\subset \ker(H_p(F^{q}X,A(n))\xrightarrow{j^*} H_p(F^{q}U,A(n))),$$ we can find $\beta\in H_p(F^qZ,A(n))$ such that $i_*(\beta)=\alpha$. Hence it suffices to show that $\beta \in H_{p}^{q,\nr}(Z,A(n))$.  Note that we have following commutative diagram by the compatibility of the residue and proper pushforward maps
\begin{equation*}
\xymatrix{H_p(F^{q}Z,A(n))\ar[r]^{i_*}\ar[d]^{\partial}& H_p(F^{q}X,A(n))\ar[d]^{\partial}\\
    \oplus_{x\in Z_{(q-1)}}H_{p-1}(x,A(n))\ar@{^(->}[r]&\oplus_{x\in X_{(q-1)}}H_{p-1}(x,A(n)),
    }
\end{equation*}
so $\partial(\beta)=0\in \oplus_{x\in Z_{(q-1)}}H_{p-1}(x,A(n))$ as $\partial(\alpha)=0\in \oplus_{x\in X_{(q-1)}}H_{p-1}(x,A(n))$ by (\ref{equ-exact sequences for homology}). That is $\beta \in H_{p}^{q,\nr}(Z,A(n))$.

Finally, the exact sequence (\ref{equ-exact sequence of F^q}) directly implies that (\ref{equ-localization sequence for refined homology}) is exact at $H_p^{q,\nr}(U,A(n))$ and at $H_{p-1}^{q-1,\nr}(Z,A(n))$. 
\end{proof}
\begin{rem}
\rm{If we change homology to cohomology, there is a canonical localization sequence
\begin{small}
  \begin{equation}\label{equ-localization sequence}
     H^{p-2c}_{q-c,\nr}(Z,A(n-c))\xrightarrow{i_\ast}  H^{p}_{q,\nr}(X,A(n))\xrightarrow{j^\ast} H^{p}_{q,\nr}(U,A(n))\xrightarrow{\delta^p_q}  H^{p+1-2c}_{q+1-c,\nr}(Z,A(n-c))\xrightarrow{i_\ast},
\end{equation}   
\end{small}

where $i\colon Z\to X$ is a closed immersion of codimension $c\coloneqq\dim(X)-\dim(Z)$ with open complement $j\colon U\hookrightarrow X$ and $d\coloneqq \dim (X)=\dim (U)$. Moreover, because we only apply the usual Gysin exact sequence in the proof of \Cref{prop-localization sequence}, (\ref{equ-localization sequence}) is valid for every twisted Borel--Moore cohomology theory on a constructible category of Noetherian schemes.
    }
\end{rem}

Moreover, there is a Mayer--Vietoris sequence for refined unramified cohomology as in \cite[Lemma 2.6]{JS03}.
\begin{prop}\label{prop-MV sequence}
Let $X$ be a Noetherian scheme and suppose $X=X_1\cup X_2$ is the union of two closed subschemes and suppose $Z\coloneqq X_1\cap X_2$. 
 The following two sequences combine and form a long exact sequence for each integer $k$
\begin{align*}
H_{p}^{p+k,\nr}(Z,A(n))\xrightarrow{((k_1)_\ast,(-k_2)_\ast)}  H_{p}^{p+k,\nr}(X_1,A(n))\oplus  H_{p}^{p+k,\nr}(X_2,A(n))\xrightarrow{(i_1)_\ast+(i_2)_\ast}  H_{p}^{p+k,\nr}(X,A(n))\\
H_{p}^{p+k,\nr}(X_1,A(n))\oplus H_{p}^{p+k,\nr}(X_2,A(n))\xrightarrow{(i_1)_\ast+(i_2)_\ast}  H_{p}^{p+k,\nr}(X,A(n))\xrightarrow{\delta}  H_{p-1}^{p+k-1,\nr}(Z,A(n)),
\end{align*}
where $k_j\colon Z\hookrightarrow X_j$ and $i_j\colon X_j\hookrightarrow X$ are the closed immersions for $j=1,2$.   
\end{prop}
\begin{proof} 
The exact sequence is induced in a standard way (the snake lemma) from the following commutative diagram of localization sequences (\ref{equ-localization sequence for refined homology})
\[
\adjustbox{scale=.9}{
\begin{tikzcd}
 H_{p}^{p+k,\nr}(Z,A(n))\ar[r,"(k_1)_\ast"]\ar[d,"(k_2)_\ast"]& H_{p}^{p+k,\nr}(X_1,A(n))\ar[r]\ar[d,"(i_1)_\ast"]&  H_{p}^{p+k,\nr}(X_1\setminus Z,A(n))\ar[r]\ar[d,"="]& H_{p-1}^{p+k-1,\nr}(Z,A(n))\ar[d,"(k_2)_\ast"]\\
 H_{p}^{p+k,\nr}(X_2,A(n))\ar[r,"(i_2)_\ast"]& H_{p}^{p+k,\nr}(X,A(n))\ar[r]&  H_{p}^{p+k,\nr}(X\setminus X_2,A(n))\ar[r]& H_{p-1}^{p+k-1,\nr}(X_2,A(n)).
\end{tikzcd}
}
\]
Moreover, the morphism $\delta$ is defined as the composition
\[
 H_p^{p+k,\nr}(X,A(n))\to  H_{p}^{p+k,\nr}(X\setminus X_2,A(n))\cong  H_{p}^{p+k,\nr}(X_1\setminus Z,A(n))\to  H_{p-1}^{p+k-1,\nr}(Z,A(n)).\qedhere
\]
\end{proof}

\begin{rem}\label{cor-ss of SNC}
 \rm{More generally, as in \cite[Proposition 3.6]{JS03} and as a formal consequence of \Cref{prop-MV sequence}, there is a spectral sequence for simplicial schemes. Let $Y$ be a Noetherian scheme and $\mathcal{Y}=(Y_1,\dots,Y_N)$ be an ordered tuple of closed subschemes with $Y=\cup Y_i$ (we do not assume any normal crossing condition, we just assume $Y$ as a simplicial scheme). Let $Y^{[s]}\coloneqq \coprod Y_{i_1,\dots,i_s}$ run through all tuples $1\leq i_1<\dots <i_s\leq N$, where $Y_{i_1,\dots,i_s}$ is the scheme-theoretic intersection of $Y_{i_1},\dots,Y_{i_s}$. For any $n$, there is a spectral sequence of homological type
\begin{align}\label{equ-ss of SNC}
		E_{s,t}^1(\mathcal{Y})= H_{t}^{t+k,\nr}(Y^{[s+1]},A(n))\Rightarrow  H^{s+t+k,\nr}_{s+t}(Y,A(n)),  
\end{align}
where 
\begin{enumerate}[(i)]
    \item $d_{s,t}^1\coloneqq \sum _{v=1} ^{s+1}(-1)^{v+1}(\delta_v)_\ast$;
    \item $(\delta_v)_\ast$ is induced by the closed immersions $Y_{i_1,\dots,i_s}\hookrightarrow Y_{i_1,\dots,\hat{i_v},\dots,i_s}$.
\end{enumerate}
}
\end{rem}
\subsection{$\mathbb{A}^1$-Homotopy Invariance}
Here we show that the refined unramified cohomology groups $H^i_{j,\nr}(-,A(n))$ satisfy the $\mathbb A^1$-homotopy invariance, provided that the ordinary twisted Borel--Moore cohomology groups do. 
\begin{define}\label{def-A1 homotopy invariance}
We say a twisted Borel--Moore cohomology $H^i(-,A(n))$ is $\mathbb{A}^1$-homotopy invariant, if it satisfies the following three conditions:
\begin{itemize}
    \item[(i)] for every flat morphism $f\colon Y\to X$ between algebraic $k$-schemes, there are functorial pullback maps $f^*\colon H^i(X,A(n))\to H^i(Y,A(n))$ compatible with localization sequences;
    \item[(ii)] for every algebraic $k$-scheme $X$, the canonical pullback $H^i(X,A(n))\to H^i(X\times \mathbb{A}^1,A(n))$ is an isomorphism for every $i$ and $n$;
    \item[(iii)] for every algebraic $k$-scheme $X$ and for every point $x\in X$, let $\eta$ be the generic point of $\mathbb{A}^1_x$ (i.e., the fiber of $x$ under the projection $X\times \mathbb{A}^1\to X$). The sequence 
     \begin{equation}\label{equ-split exact for homotopy invariant}
         0\to H^i(x,A(n))\xrightarrow{r}H^{i}(\eta,A(n))\xrightarrow{d} \bigoplus_{p\in {{\mathbb{A}^1_{x}}}^{(1)}}H^{i-1}(p,A(n-1))\to 0
     \end{equation}
     is exact, where $p$ runs through all $1$-codimensional points of $\mathbb{A}^1_x$, $r$ is induced by the flat pull-back $X\times \mathbb{A}^1\to X$ and $d$ is induced by the connecting map.
\end{itemize}
\end{define}
\begin{rem}
\rm{Let us elaborate on (iii) from \Cref{def-A1 homotopy invariance} here. This appears in \cite[Proposition 2.2]{Ros96} and \cite[Proposition 4.1.4]{colliot1995birational} (in their cases, (\ref{equ-split exact for homotopy invariant}) is split exact) naturally, and it is the essence of the proof of the $\mathbb{A}^1$-homotopy invariance of $E^{p,q}_2$, which the reader will see from the proof of \Cref{lem-A1 homopoty invariant}.
In fact, applying \Cref{prop-Gysin sequence} for every non-empty open subset $U\subset X$, there is a long exact sequence
\begin{small}
\begin{equation*}
   \to H^i(U\times \mathbb{A}^1,A(n))\to H^i(F_0(U\times \mathbb{A}^1),A(n))\xrightarrow{\partial}\bigoplus_{p\in (U\times \mathbb{A}^1)^{(1)}}H^{i-1}(p,A(n-1))\xrightarrow{\iota_*} H^{i+1}(U\times \mathbb{A}^1,A(n))\to.
\end{equation*}
\end{small}

Taking direct limits over $U$, and noting that $H^i(U\times \mathbb{A}^1,A(n))\cong H^i(U,A(n))$ by flat pull-back (ii), then we have the following exact sequence
\begin{equation*}
  \cdots\to H^i(x,A(n))\xrightarrow{r}H^{i}(\eta,A(n))\xrightarrow{d} \bigoplus_{p\in {{\mathbb{A}^1_{x}}}^{(1)}}H^{i-1}(p,A(n-1))\to  H^{i+1}(x,A(n))\to\cdots.
\end{equation*}
Therefore, the condition (iii) is equivalent to the injectivity of $r$ or to the surjectivity of $d$. 
}
\end{rem}
\begin{prop}
 All examples in \Cref{section-2.1}, (i.e., (a)-(d) in \Cref{section-2.1}) are twisted Borel--Moore cohomology theories with $\mathbb{A}^1$-homotopy invariance.   
\end{prop}
\begin{proof}
It suffices to show the morphism $r$ in (\ref{equ-split exact for homotopy invariant}) is injective.  For the theory (a), this is already given by the theory of cycle modules (see \cite[Remark 1.11, Proposition 2.2]{Ros96}). More specifically,  observing that pro-\'etale cohomology does not change by purely inseparable extensions (see \cite[Lemma 5.4.2]{BS15}), we may assume that $k$ is perfect so that there exists an open non-empty smooth subset $U\subset \overline{\{x\}}$. By \cite[Lemma 6.5]{schreieder2021refined} we have isomorphisms
\begin{equation*}
H^i(U,\mu_m^{\otimes n})\cong H^i_{\text{cont}}(U,(\mu_m^{\otimes n},\id)_s)\cong H^i_{\e t}(U,\mu_m^{\otimes n}),    
\end{equation*}
where the middle one is the \emph{continuous \'etale cohomology} with the constant inverse system $(\mu_m^{\otimes n},\id)_s$. Taking direct limits over $U$, we get $H^i(x,\mu_m^{\otimes n})\cong  H^i_{\e t}(x,\mu_m^{\otimes n})$ for every $i$, $n$ and $m$. For a fixed $\kappa(x)$-rational point $Q$ on $\mathbb{A}^1_x$, let $\partial$ be the residue map $H^*(\eta,A(\dagger))\to H^{*-1}(Q,A(\dagger-1))$ and $\pi$ be a prime of the valuation on $\kappa(\eta)$ induced by $Q$. Let $\epsilon\colon\kappa(\eta)^*\to H^1(\eta,\Z_{\ell}(1))$ as in \cite[Definition 4.4, P6]{schreieder2021refined}. Then \cite[Definition~1.1,~R3d]{Ros96} implies that $\partial(\epsilon(\pi)\cup r(\alpha))=\alpha$ by specializing for all $\alpha\in H^i(x,A(n))$, so $r$ is injective.  Similarly for theory (d), consider the map $\epsilon\colon\kappa(\eta)^*\to H^1(\eta,\Z(1))$ mentioned in the proof of \cite[Proposition 6.9]{schreieder2020infinite}, we still have $\partial(\epsilon(\pi)\cup r(\alpha))=\alpha$ for $\alpha\in H^i(x,A(n))$ and then $r$ is injective as well. In particular, the exact sequence (\ref{equ-split exact for homotopy invariant}) is split in these cases.

For the case (b), again by \cite[Lemma 6.5]{schreieder2021refined}, there is an isomorphism $$H^i(x,\Z_{\ell}(n))\cong H^{i}_{\text{cont}}(x,(\mu_{\ell^r}^{\otimes n})_{r})\eqqcolon H^i_{\text{cont}}(x,\Z_{\ell}(n)).$$ Then the following exact sequence 
  \begin{equation*}
    0\to R^1\varprojlim_r H^{i-1}_{\e t}(x,\Z/\ell^r(n))\to H^i_{\text{cont}}(x,\Z_{\ell}(n))\to \varprojlim_r H^i_{\e t}(x,\Z/\ell^r(n))\to 0
  \end{equation*}
  and the splitting of (\ref{equ-split exact for homotopy invariant}) for (a) gives the injectivity of $r$ in case (b). Finally, (c) is $\A^1$-homotopy invariant because (a) and (b) now are.\qedhere
\end{proof}
\begin{lem}[cf.~{\cite[Proposition 8.6]{Ros96}}]\label{lem-A1 homopoty invariant}
    Let $X$ be an algebraic $k$-scheme and let $Y\coloneqq X\times \mathbb{A}^1$ with the projection $\pi:Y\to X$. If the twisted Borel--Moore cohomology $H^i(-,A(\dagger))$ is $\mathbb{A}^1$-homotopy invariant, then the pull-back $E_2^{i,j+i}(X,A(n+i))\to E_2^{i,j+i}(Y,A(n+i))$ is an isomorphism for all $i,j$ and $n$, where $E_2^{i,j+i}(*,A(n+i))$ is defined in \Cref{prop-exact sequence of D2}.
\end{lem}
\begin{proof}
For every $y\in Y^{(i)}$, let $u\coloneqq \pi(y)$, and then $u\in X^{(i)}$ or $u\in X^{(i-1)}$. If $u\in X^{(i)}$, $y$ is a $0$-codimensional point in $\mathbb{A}^1_u\coloneqq \text{Spec} (\kappa(u)[t])$ over $u$, so $y$ is exactly the generic point of $\mathbb{A}^1_u$, denoted by $\eta(u)$; if $u\in X^{(i-1)}$, $y$ is a $1$-codimensional point in $\mathbb{A}^1_u$ over $u$. Therefore, we have the following equations
\begin{eqnarray*}
E_1^{i,j+i}(Y,A(n+i))\coloneqq\bigoplus_{y\in Y^{(i)}}H^{j}(y,A(n))=\bigoplus_{\pi(y)\in X^{(i)}}H^{j}(y,A(n))\oplus \bigoplus_{\pi(y)\in X^{(i-1)}}H^{j}(y,A(n))\notag\\
= \bigoplus_{u\in X^{(i)}}H^{j}(\eta(u),A(n))\oplus\bigoplus_{u\in X^{(i-1)}}\bigoplus_{y\in {\mathbb{A}^1_{u}}^{(1)}}H^{j}(y,A(n)).
\end{eqnarray*}
 So $E_2^{i,j+i}(Y,A(n+i))$ can be computed by a double complex (up to a sign of the differential) with only two non-trivial rows whose ${'D_0}^{p,q}$ looks like
\begin{equation*}
\adjustbox{scale=.8}{
    \xymatrix{ \bigoplus\limits_{u\in X^{(i-1)}}H^{j+1}(\eta(u),A(n+1))\ar[r]^{\partial}\ar@{->>}[d]^{\partial}&\bigoplus\limits_{u\in X^{(i)}}H^{j}(\eta(u),A(n))\ar[r]^{\partial}\ar@{->>}[d]^{\partial}&\bigoplus\limits_{u\in X^{(i+1)}}H^{j-1}(\eta(u),A(n-1))\ar@{->>}[d]^{\partial}\\
\bigoplus\limits_{u\in X^{(i-1)}}\bigoplus\limits_{y\in {\mathbb{A}^1_{u}}^{(1)}}H^{j}(y,A(n))\ar[r]^{\partial}&\bigoplus\limits_{u\in X^{(i)}}\bigoplus\limits_{y\in {\mathbb{A}^1_{u}}^{(1)}}H^{j-1}(y,A(n-1))\ar[r]^{\partial}&\bigoplus\limits_{u\in X^{(i+1)}}\bigoplus\limits_{y\in {\mathbb{A}^1_{u}}^{(1)}}H^{j-2}(y,A(n-2)),
    }}
\end{equation*}
where all vertical arrows $\partial$ are induced by $d$ in (\ref{equ-split exact for homotopy invariant}) and are surjective. More precisely, we have 
\begin{equation*}
  'D_0^{p,q}=\left\{\begin{array}{ll}
      \bigoplus\limits_{u\in X^{(p)}}H^{i+j-p}(\eta(u),A(n+i-p)) & \text{if } q=0 \\
      \bigoplus\limits_{u\in X^{(p)}}\bigoplus\limits_{y\in {\mathbb{A}^1_{u}}^{(1)}}H^{i+j-1-p}(y,A(n+i-p-1)) & \text{if } q=1\\
      0 & \text{else.}
  \end{array} \right.
\end{equation*}
This double complex is compatible with the flat pull-back $\pi^*$.
Let $'d_0^{p,q}:{'D_0^{p,q}}\to {'D_0^{p,q+1}}$ be the vertical arrows. Applying (\ref{equ-split exact for homotopy invariant}), we get that
\begin{equation*}
  'E_1^{p,q}=\left\{\begin{array}{ll}
      \bigoplus_{u\in X^{(p)}}H^{i+j-p}(u,A(n+i-p)) & \text{if } q=0 \\
      0 & \text{else.}
  \end{array} \right.
\end{equation*}
Moreover, $'d_1^{p,0}:{'E_1^{p,0}\to {'E_1^{p+1,0}}}$ is $d_1^{p,i+j}:E_1^{p,i+j}(X,A(n+i))\to E_1^{p+1,i+j}(X,A(n+i))$ by compatibility of localization sequences and flat pull-backs. So taking cohomology we see the above spectral sequence $\{'E_r^{p,q},'d_r\}_r$ degenerates at page $2$ and for all $i$ and $j$, we have $$E_2^{i,i+j}(X,A(n+i))\xrightarrow[\sim]{\pi^*} {'E_2^{i,0}}={'E_{\infty}^{i,0}}=H^i('D_0^{\bullet,\bullet})=E_2^{i,i+j}(Y,A(n+i)).$$
In other words, the map $\pi^*\colon E_2^{i,j+i}(X,A(n+i))\to E_2^{i,j+i}(Y,A(n+i))$ is an isomorphism for all $i$ and $j$.\qedhere

\end{proof}
\begin{prop}\label{prop-A1 homotopy invariance}
Let $X$ be an algebraic $k$-scheme and let $\pi\colon Y\to X$ be an affine bundle. If the twisted Borel--Moore cohomology $H^i(-,A(n))$ is $\mathbb{A}^1$-homotopy invariant, then $H^i_{j,\nr}(-,A(n))$ is as well. That is, the pull-back $\pi^\ast\colon  H^i_{j,\nr}(X,A(n))\to  H^i_{j,\nr}(Y,A(n))$ is an isomorphism for all $i$ and $j$.
\end{prop}
\begin{proof}
We first reduce it to a simple case. By \Cref{prop-localization sequence}, applying five-lemma and arguments involving dimensional induction, we may assume that $Y$ is a trivial vector bundle over $X$, i.e., $Y\cong X\times \mathbb{A}^r$ for some $r$. Furthermore, we may assume that $Y=X\times \mathbb{A}^1$ by the functoriality of flat pullbacks.

Let us do induction on $j$. For $j$ large (i.e., $j\geq \dim(X)+1$), it follows from the equalities $H^i_{j,\nr}(X,A(n))=H^i(X,A(n))$, $H^i_{j,\nr}(Y,A(n))=H^i(Y,A(n))$ and the homotopy invariance of $H^i(-,A(n))$. Now suppose that for any $p\geq j$ and all $i$, the $\pi^\ast\colon  H^i_{p,\nr}(X,A(n))\to  H^i_{p,\nr}(Y,A(n))$ is an isomorphism. We show that $\pi^\ast\colon  H^i_{j-1,\nr}(X,A(n))\to  H^i_{j-1,\nr}(Y,A(n))$ is an isomorphism.

The pull-back map $\pi^\ast$ induces a morphism of long exact sequences from \Cref{prop-exact sequence of D2}
\[
\adjustbox{scale=.8}{
\begin{tikzcd}
E_2^{j,j-i}(X,A(n))\ar[r]\ar[d,"\pi^\ast"]& H^i_{j,\nr}(X,A(n))\ar[r]\ar[d,"\pi^\ast"]&  H^{i}_{j-1,\nr}(X,A(n))\ar[r]\ar[d,"\pi^\ast"]&E_2^{j+1,i-j}(X,A(n))\ar[r]\ar[d,"\pi^\ast"]& H^{i+1}_{j+1,\nr}(X,A(n))\ar[d,"\pi^\ast"]\\
E_2^{j,j-i}(Y,A(n))\ar[r]& H^i_{j,\nr}(Y,A(n))\ar[r]&  H^{i}_{j-1,\nr}(Y,A(n))\ar[r]&E_2^{j+1,i-j}(Y,A(n))\ar[r]& H^{i+1}_{j+1,\nr}(Y,A(n))\nospacepunct{.}
\end{tikzcd}
}
\]
As $E_2^{p,q}$ has the $\A^1$-homotopy invariance by \Cref{lem-A1 homopoty invariant}, we conclude by induction and the five-lemma that the middle arrow is an isomorphism as well. 
\end{proof}
\begin{rem}\label{rem-homology for An}
\rm{A direct corollary of \Cref{prop-A1 homotopy invariance} is that 
\begin{equation*}
  H^i_{j,\nr}(\mathbb{A}^m_k,A(n))\cong H^i_{j,\nr}(\mathrm{Spec}(k),A(n))=\left\{\begin{array}{ll}
      H^i(\mathrm{Spec}(k),A(n)) & \text{if } j\geq 0 \\
      0 & \text{if } j<0
  \end{array} \right.
\end{equation*}
by \Cref{lem-canonical morphism is an isomorpshim},  while a direct computation of $H^i_{j,\nr}(\mathbb{A}^m_k,A(n))$ appears to be a rather difficult task, because $H^i(F_j\mathbb{A}^m_k,A(n))$ is typically a large non-finitely generated group.
}
\end{rem}

\subsection{Bockstein Sequence}\label{section-Bockstein sequence}
For twisted Borel--Moore cohomology theories, there exists a Bockstein long exact sequence
\begin{equation}\label{equ-original Bockstein sequence}
\cdots\to H^i(X,\Z_\ell(n))\xrightarrow{\times\ell^r} H^i(X,\Z_\ell(n))\to H^i(X,\mu_{\ell^r}^{\otimes n})\xrightarrow{\delta} H^{i+1}(X,\Z_\ell(n))\xrightarrow{\times\ell^r}\cdots,    
\end{equation}
which is functorial with respect to open restriction and proper pushforward, \cite[Definition 4.4(P5)]{schreieder2021refined}. In this subsection we show there is a natural generalisation to refined unramified cohomology if this twisted Borel--Moore cohomology satisfies a Bloch--Kato type conjecture. Note that if $k=\C$, then we can replace $\Z_\ell$ and $\mu_{\ell^r}$ by $\Z$ and $\Z/n\Z$ for any $n\in \Z$, respectively. 
\begin{define}[{cf.~\cite[Remark 5.14]{schreieder2021refined}}]\label{def-BK type}
We say a twisted Borel--Moore cohomology $H^i(-,A(n))$ satisfies Bloch--Kato $\ell$-type conjecture, if for any algebraic $k$-scheme $X$ and $x\in X$,  there are surjections $(\kappa(x)^*)^{\otimes i}\twoheadrightarrow H^i(x,\mu_{\ell^r}^{\otimes i})$ that factor through $H^i(x,\Z_{\ell}(i))$ for all $r$ and $i$.
\end{define}
\begin{rem}\label{rem-torsion free}
\rm{By Bloch's observation, if $H^i(-,A(n))$ satisfies Bloch--Kato $\ell$-type conjecture, then $H^{i+1}(x,\Z_{\ell}(i))$ is torsion free for all $i$ (see \cite{bloch2010lectures}). Moreover, all examples in \Cref{section-2.1}, (i.e., (a)-(d) in \Cref{section-2.1}) are Borel--Moore cohomology theories satisfying Bloch--Kato $\ell$-type conjecture for all primes $\ell\in k^*$ by the \emph{Bloch--Kato conjecture} (see \Cref{thm-Bloch Kato conjecture}) and Poincar\'e duality. 
}
\end{rem}

\begin{lem}\label{lem-construt bockstein sequence}
Suppose the twisted Borel--Moore cohomology $H^i(-,A(n))$ satisfies Bloch--Kato $\ell$-type conjecture. Then for any algebraic $k$-scheme $X$, there is a functorial morphism 
\[
\Delta_j^i\colon H^i_{j,\nr}(X,\mu_{\ell^r}^{\otimes (i-j-1)})\to H^{i+1}_{j+1,\nr}(X,\Z_\ell(i-j-1)),
\]
compatible with the $F_j$ filtrations on refined unramified cohomology, which is the original Bockstein map $\delta$ for either $j\geq \lceil i/2\rceil$ or $j\geq \dim(X)$.
\end{lem}
\begin{proof}
Let $\alpha\in H^i_{j,\nr}(X,\mu_{\ell^r}^{\otimes(i-j-1)})\subset H^i(F_{j}X,\mu_{\ell^r}^{\otimes (i-j-1)})$, by definition, $\alpha$ can be lifted to an element $\alpha'\in H^i(F_{j+1}X,\mu_{\ell^r}^{\otimes (i-j-1)})$, and the Bockstein connecting map gives $\delta(\alpha')$ in  $H^{i+1}(F_{j+1}X,\Z_{\ell}(i-j-1))$. Then $\partial(\delta(\alpha'))=\delta(\partial(\alpha'))\in\bigoplus_{x\in X^{(j+2)}} H^{i-2j-2}(x,\Z_\ell(i-2j-3))$ is torsion, hence zero since $ H^\ast(x,\Z_{\ell}(\ast-1))$ is torsion-free by \Cref{rem-torsion free}. So we conclude that $\delta(\alpha')\in  H^{i+1}_{j+1,\nr}(X,\Z_{\ell}(i-j-1))$ and we set $\Delta^i_j(\alpha)\coloneqq \delta(\alpha')$.

We are left to show that this definition does not depend on the lift $\alpha'$ of $\alpha$. So suppose $\alpha''\in  H^i(F_{j+1}X,\mu_{\ell^r}^{\otimes (i-j-1)})$ is another lift, then $\alpha'-\alpha''$ comes from some $\beta\in \bigoplus_{x\in X^{(j+1)}} H^{i-2j-2}(x,\mu_{\ell^r}^{\otimes (i-2j-2)})$ by \Cref{prop-Gysin sequence}. Consider the following commutative diagram
\[
\begin{tikzcd}
\bigoplus\limits_{x\in X^{(j+1)}} H^{i-2j-2}(x,\mu_{\ell^r}^{\otimes (i-2j-2)})\ar[r,"\iota_\ast"]\ar[d,"\delta=0"]& H^i(F_{j+1}X,\mu_{\ell^r}^{\otimes (i-j-1)})\ar[d,"\delta"]\\
\bigoplus\limits_{x\in X^{(j+1)}} H^{i-2j-1}(x,\Z_{\ell}(i-2j-2))\ar[r,"\iota_\ast"]& H^{i+1}(F_{j+1}X,\Z_{\ell}(i-j-1))
\end{tikzcd}
\]
where the left arrow is zero as $H^{i-2j-1}(x,\Z_{\ell}(i-2j-2))$ is torsion-free by \Cref{rem-torsion free} again. So $\delta(\alpha'-\alpha'')=\delta(\iota_\ast(\beta))=\iota_\ast(\delta(\beta))=0$, and $\Delta_j^i(\alpha)=\delta(\alpha')=\delta(\alpha'')$ is well-defined. 

The functoriality of $\Delta^i_j$ follows directly from the same functoriality properties of the original Bockstein sequence (\ref{equ-original Bockstein sequence}). Finally by \Cref{cor-degeneration of refined unramified cohomology} and \Cref{def-refined unramified cohomology}, the restriction morphisms $H^i(X,\mu_{\ell^r}^{\otimes n})\to H^i_{j,\nr}(X,\mu_{\ell^r}^{\otimes n})$ and $H^i(X,\Z_{\ell}(n))\to H^i_{j,\nr}(X,\Z_{\ell}(n))$ are both isomorphisms, for either $j\geq \lceil i/2\rceil$ or $j\geq \dim(X)$. We finish our proof by noting that $j\geq \lceil i/2\rceil$ implies $j+1\geq \lceil (i+1)/2\rceil$ and the following commutative diagram making use of the functoriality properties of the original Bockstein sequence
\begin{equation*}
\xymatrix{H^i(F_{j+1}X,\mu_{\ell^r}^{\otimes (i-j-1)})\ar[r]^{\delta}\ar[d]^{\rm restr.}&H^{i+1}(F_{j+1}X,\Z_{\ell}(i-j-1))\ar[d]^{\rm restr.}\\
H^i(F_{j}X,\mu_{\ell^r}^{\otimes (i-j-1)})\ar[r]^{\delta}&H^{i+1}(F_{j}X,\Z_{\ell}(i-j-1)).
}
\end{equation*}
\end{proof}
\begin{rem}
\rm{In this proof, to define $\Delta_j^i$, the torsion-freeness of $H^{i-2j-1}(x,\Z_{\ell}(i-2j-2))$ is enough. By \cite[Lemma 5.13]{schreieder2021refined}, we therefore can define $\Delta_j^i$ for all $i,j$ with $i-2j\leq 3$ for every twisted Borel-Moore cohomology theory.  
 }
\end{rem}

\begin{prop}\label{prop-short Bockstein sequence}
Let $X$ be an algebraic $k$-scheme.  Suppose the twisted Borel--Moore cohomology $H^i(-,A(n))$ satisfies Bloch--Kato $\ell$-type conjecture. There is a canonical long exact sequence for refined unramified cohomology for all $i,j\geq 0$
\begin{small}
\begin{equation*}
  \to H^i_{j,\nr}(X,\Z_\ell(i-j-1))\xrightarrow{\times \ell^r}H^i_{j,\nr}(X,\Z_\ell(i-j-1))\to  H^i_{j,\nr}(X,\mu_{\ell^r}^{\otimes (i-j-1)})\xrightarrow{\Delta^i_j} H^{i+1}_{j+1,\nr}(X,\Z_\ell(i-j-1))\xrightarrow{\times \ell^r}.      
\end{equation*}
\end{small}
\end{prop}
\begin{proof}
First the exactness at $ H^i_{j,\nr}(X,\mu_{\ell^r}^{\otimes (i-j-1)})$ follows directly from the exactness of $$ H^i(F_{j+1}X,\Z_{\ell}(i-j-1))\to  H^i(F_{j+1}X,\mu_{\ell^r}^{\otimes (i-j-1)})\xrightarrow{\delta} H^{i+1}(F_{j+1}X,\Z_{\ell}(i-j-1)).$$ 

Then we show that $$ H^i_{j,\nr}(X,\mu_{\ell^r}^{\otimes (i-j-1)})\xrightarrow{\Delta^i_j} H^{i+1}_{j+1,\nr}(X,\Z_{\ell}(i-j-1))\xrightarrow{\times \ell^r} H^{i+1}_{j+1,\nr}(X,\Z_{\ell}(i-j-1))$$ is exact. Its composition is zero as $$ H^i(F_{j+1}X,\mu_{\ell^r}^{\otimes (i-j-1)})\to  H^{i+1}(F_{j+1}X,\Z_{\ell}(i-j-1))\to  H^{i+1}(F_{j+1}X,\Z_{\ell}(i-j-1))$$ is exact. Now let $$\alpha\in \ker( H^{i+1}_{j+1,\nr}(X,\Z_{\ell}(i-j-1))\xrightarrow{\times \ell^r} H^{i+1}_{j+1,\nr}(X,\Z_{\ell}(i-j-1))\subseteq  H^{i+1}(F_{j+1}X,\Z_{\ell}(i-j-1));$$ then by the exactness of the aforementioned sequence, there exists a $\beta'\in  H^i(F_{j+1}X,\mu_{\ell^r}^{\otimes (i-j-1)})$ with $\delta(\beta')=\alpha$. This $\beta'$ defines an element $\beta\in H_{j,\nr}^i(X,\mu_{\ell^r}^{\otimes (i-j-1)})$ and by definition of $\Delta^i_j$, we see that $\Delta^i_j(\beta)=\delta(\beta')=\alpha$, as wished.

Finally, we show the exactness at $H^i_{j,\nr}(X,\Z_{\ell}(i-j-1))$.  Note that this sequence is induced by the sequence 
\[
\cdots\to H^i(F_jX,\Z_\ell(i-j-1))\xrightarrow{\times \ell^r} H^i(F_jX,\Z_\ell(i-j-1))\to H^i(F_jX,\mu_{\ell^r}^{\otimes (i-j-1)})\xrightarrow{\delta} \cdots,
\]
hence the composition is zero. Let $$\alpha\in \ker( H^i_{j,\nr}(X,\Z_{\ell}(i-j-1))\to  H^i_{j,\nr}(X,\mu_{\ell^r}^{\otimes (i-j-1)}))\subset H^i(F_jX,\Z_{\ell}(i-j-1)).$$ Then $\alpha$ maps to zero under the morphism $H^i(F_jX,\Z_\ell(i-j-1))\to H^i(F_jX,\mu_{\ell^r}^{\otimes (i-j-1)})$, so we can find $\beta\in H^i(F_jX,\Z_\ell(i-j-1))$ such that $\ell^r\cdot\beta=\alpha\in H^i(F_jX,\Z_\ell(i-j-1))$ using the aforementioned exact sequence. This implies that $$\ell^r\cdot\partial(\beta)=\partial(\ell^r\cdot\beta)=\partial(\alpha)=0\in \bigoplus_{x\in X^{(j+1)}}H^{i-2j-1}(x,\Z_{\ell}(i-2j-2)),$$
hence $\partial(\beta)=0$ by the torsion-freeness of $H^{i-2j-1}(x,\Z_{\ell}(i-2j-2))$, hence $\beta$ belongs to $H^i_{j,\nr}(X,\Z_{\ell}(i-j-1))$ by \cite[Lemma 5.8]{schreieder2021refined}. So we get 
\begin{small}
\begin{equation*}
    \ker( H^i_{j,\nr}(X,\Z_{\ell}(i-j-1))\to  H^i_{j,\nr}(X,\mu_{\ell^r}^{\otimes (i-j-1)}))\subset\im(H^i_{j,\nr}(X,\Z_{\ell}(i-j-1))\to  H^i_{j,\nr}(X,\mu_{\ell^r}^{\otimes( i-j-1)})).\qedhere
\end{equation*}    
\end{small}
\end{proof}

\begin{cor}\label{prop-Bockstein sequence}
Let $X$ be an algebraic scheme defined over a field containing all $\ell$-power roots of unity. Suppose that the twisted Borel--Moore cohomology $H^i(-,A(n))$ are one of examples in \Cref{section-2.1}. There is a long exact Bockstein sequence for refined unramified cohomology
\[
\cdots\to  H^i_{j,\nr}(X,\Z_\ell(n))\xrightarrow{\times \ell^r}  H^i_{j,\nr}(X,\Z_\ell(n))\to  H^i_{j,\nr}(X,\mu_{\ell^r}^{\otimes n})\xrightarrow{\Delta^i_j} H^{i+1}_{j+1,\nr}(X,\Z_\ell(n))\xrightarrow{\times \ell^r}\cdots.   
\]
\end{cor}
\begin{proof}
The morphism $H^i_{j,\nr}(X,\mu_{\ell^r}^{\otimes n})\xrightarrow{\Delta^i_j} H^{i+1}_{j+1,\nr}(X,\Z_\ell(n))$ is defined by \Cref{lem-construt bockstein sequence} after choosing isomorphisms of sheaves $\mu_{\ell^r}^{\otimes n}\cong \mu_{\ell^r}^{\otimes (i-j-1)}$ and $\Hat{\Z}_\ell(n)\cong \Hat{\Z}_\ell(i-j-1)$. We finish the proof by applying \Cref{prop-short Bockstein sequence}.
\end{proof}
\begin{rem}
{\rm The proof of \Cref{prop-short Bockstein sequence} is actually a proof of \cite[Theorem 7.7]{schreieder2021refined} in disguise. Namely, we see directly that $Z^i(X)[\ell^r]\cong H_{i-1,\nr}^{2i}(X,\Z_\ell(i))[\ell^r]$ for an invertible prime $\ell$ and the Bockstein sequence gives that this is isomorphic to $\frac{H_{i-2,\nr}^{2i-1}(X,\mu_{\ell^r}^{\otimes i})}{H_{i-2,\nr}^{2i-1}(X,\Z_\ell(i))}$ via $\Delta^{2i-1}_{i-2}$. Similarly, now taking the direct limit over $r$, we have $$Z^i(X)[\ell^\infty]\cong H_{i-1,\nr}^{2i}(X,\Z_\ell(i))[\ell^\infty]\cong\frac{H_{i-2,\nr}^{2i-1}(X,\Q_\ell/\Z_\ell(i))}{H_{i-2,\nr}^{2i-1}(X,\Q_\ell(i))}.$$}
\end{rem}

Taking the direct limit over $r$ in \Cref{prop-Bockstein sequence}, we conclude the following.
\begin{cor}\label{cor-Bockstein sequence Q/Z}
Let $X$ be an algebraic scheme defined over a field containing all $\ell$-power roots of unity. Suppose $H^i(X,A(n))$ is one of the examples in \Cref{section-2.1}. There is a natural long exact sequence
\[
\cdots\to  H^i_{j,\nr}(X,\Z_\ell(n))\xrightarrow{}  H^i_{j,\nr}(X,\mathbb{Q}_\ell(n))\to  H^i_{j,\nr}(X,\mathbb{Q}_\ell/\Z_\ell(n))\xrightarrow{} H^{i+1}_{j+1,\nr}(X,\Z_\ell(n))\to\cdots.   
\]
\end{cor}
\subsection{The Refined Unramified Cohomology of $X \times \mathbb{P}^r$}
In this section, we compute the refined unramified cohomology groups of $X\times\P^r$ for every algebraic $k$-scheme $X$ as a corollary of \Cref{prop-localization sequence} and \Cref{prop-A1 homotopy invariance}. 

Let $H\coloneqq X\times\P^{r-1}$ be a general (relative) hyperplane section and write $h\colon H\hookrightarrow X\times \P^r$ for the inclusion. Moreover, set $\pi\colon X\times\P^r\to X$ for the projection and $\pi'\colon H\to X$ for its restriction. More generally, we write $H^k=X\times\P^{r-k}$ for the intersection of $k$ general hyperplane sections satisfying $H^{r}\subset H^{r-1}\subset\cdots\subset H^2\subset H^1$ and $(\pi^k)'\colon H^k\to X$ and $h^k\colon H^k\hookrightarrow X\times \P^r$ for its (flat) projection and inclusion respectively.

\begin{define}
Suppose the twisted Borel--Moore cohomology theory $H^i(-,A(n)$ is $\mathbb{A}^1$-homotopy invariant. We define $P^k\colon  H^p_{q,\nr}(X,A(n))\to  H^{p+2k}_{q+k,\nr}(X\times\P^r,A(n+k))$ as the composition 
\[
H_{q,\nr}^p(X,A(n))\xrightarrow{((\pi^k)')^\ast} H_{q,\nr}^p(H^k,A(n))\xrightarrow{h^k_\ast}H_{q+n,\nr}^{p+2k}(X\times\P^r,A(n+k)).
\]
We also set $P^0\coloneqq \pi^\ast$ and set $P\coloneqq P^1$ for simplicity.
\end{define}

\begin{rem}
{\rm We warn the reader that for an abstract twisted Borel--Moore cohomology theory with $\mathbb{A}^1$-homotopy invariance, the morphism $P^k$ may depend on the choice of $H^k$. Even for $r=1$, let $Q$ and $Q'$ be two different points of $\mathbb{P}^1$, we may have $$P_{Q}\neq P_{Q'}:H^p_{q,\nr}(Spec(k),A(n))\to H^{p+2}_{q+1,\nr}(\mathbb{P}^1,A(n+2)).$$
However, for all twisted Borel--Moore cohomology theories mentioned in \Cref{section-2.1}, the $P^k$ is independent of the choice of $H^k$ and we do not need the inclusion $H^{r}\subset\cdots\subset  H^1$, because they have cup products and cycle class maps $\CH_i(*)\to H_{2i}(*,A(2i))$. Therefore these Borel--Moore cohomology theories have bi-additive Chow groups actions and $[H^k]=[(H^k)']\in \CH^k(X\times \mathbb{P}^r)$, where $P^k$ is exactly the cycle action $[H^k]_*$. Moreover, if $X$ is smooth projective, then we note that this map is the same as the composition 
\[
H_{q,\nr}^p(X,A(n))\xrightarrow{\pi^\ast}H_{q,\nr}^p(X\times\P^r,A(n))\xrightarrow{[H^k]_\ast}H_{q+k,\nr}^{p+2k}(X\times\P^r,A(n+k)),
\]
where $[H^k]_\ast$ is the action of correspondence by $\Delta_\ast([H^k])\in \CH^{k+d_X+r}(X\times\P^r\times X\times\P^r)$ as defined in \cite[Corollary 6.8]{schreieder2022moving}.}
\end{rem}
\begin{prop}[{cf.~\cite[Corollary 6.11]{schreieder2022moving}}]\label{computation trivial proj bundle}
Let $X$ be an algebraic $k$-scheme and $H^i(-,A(n))$ be a twisted Borel--Moore cohomology theory with $\mathbb{A}^1$-homotopy invariance. The map
\[
\bigoplus_{k=0}^r H_{q-k,\nr}^{p-2k}(X,A(n-k))\xrightarrow{\sim}H_{q,\nr}^p(X\times\P^r,A(n))\colon(\alpha_0,\dots,\alpha_r)\mapsto\sum_{k=0}^rP^k(\alpha_k)
\]
is an isomorphism.
\end{prop}
\begin{proof}
We do induction on $r$. For $r=0$, this is trivial, so assume the statement holds for $X\times\P^r$, and then we show that it holds for $X\times \P^{r+1}$. 

Consider the (trivial relative) hyperplane section $i\colon X\times \P^r\hookrightarrow X\times \P^{r+1}$, with open complement $j\colon X\times \mathbb{A}^{r+1}\hookrightarrow X\times \P^{r+1}$. By \Cref{prop-localization sequence}, we have the long exact sequence
\[
H^{p-2}_{q-1,\nr}(X\times \P^{r},A(n-1))\xrightarrow{i_\ast}  H^{p}_{q,\nr}(X\times \P^{r+1},A(n))\xrightarrow{j^\ast} H^{p}_{q,\nr}(X\times \mathbb{A}^{r+1},A(n))\xrightarrow{\delta^p_q}  H^{p-1}_{q,\nr}(X\times \P^{r},A(n-1)).
\]
By \Cref{prop-A1 homotopy invariance}, $(\pi\circ j)^\ast\colon  H^p_{q,\nr}(X,A(n))\to  H^{p}_{q,\nr}(X\times \mathbb{A}^{r+1},A(n))$ is an isomorphism, hence $j^\ast$ is surjective. This gives the short exact sequence
\begin{equation*}
    0\to  H^{p-2}_{q-1,\nr}(X\times \P^r,A(n-1))\xrightarrow{i_\ast}  H^{p}_{q,\nr}(X\times \P^{r+1},A(n))\to  H^{p}_{q,\nr}(X,A(n))\to 0,
\end{equation*}
which splits via the section map $\pi^\ast\colon  H^p_{q,\nr}(X,A(n))\to  H^p_{q,\nr}(X\times\mathbb{P}^{r+1},A(n))$. Thus 
\begin{equation*}
 \begin{split}
     H^p_{q,\nr}(X\times \P^{r+1},A(n))&\xleftarrow{\sim}  H^p_{q,\nr}(X,A(n))\oplus  H^{p-2}_{q-1,\nr}(X\times \P^r,A(n-1))\\
     &\xleftarrow{\sim}  H^p_{q,\nr}(X,A(n))\oplus  H^{p-2}_{q-1,\nr}(X,A(n-1))\oplus \cdots \oplus  H^{p-2r-2}_{q-r-1,\nr}(X,A(n-r-1)),
 \end{split}   
\end{equation*}
where the isomorphism is precisely given by $(\pi^\ast,P,\dots,P^{r+1})$.
\end{proof}

\begin{prop}\label{prop-G_m bundle}
Let $X$ be an algebraic $k$-scheme and $H^i(-,A(n))$ be a twisted Borel--Moore cohomology theory with $\mathbb{A}^1$-homotopy invariance. There is a canonical exact sequence
\begin{equation}\label{equ-G_m bundle}
    0\to H^p_{q,\nr}(X,A(n))\xrightarrow{\pi^*}H^p_{q,\nr}(X\times \mathbb{G}_m,A(n))\xrightarrow{\delta^p_q} H^{p-1}_{q,\nr}(X,A(n-1))\to 0,
\end{equation}
 where $\pi\colon X\times \mathbb{G}_m\to X$ is the flat projection and $\delta^p_q$ is induced by the closed immersion $X\times 0\hookrightarrow X\times \mathbb{A}^1$ by \Cref{lem-construct delta}.
\end{prop}
\begin{proof}
Let us consider the following commutative diagram of localization exact sequences (we leave out the coefficients $A$ for simplicity)
\begin{small}
\begin{equation*}
\xymatrix{H^{p-2}_{q-1,\nr}(X\times \{0\},n-1)\ar@{^(->}[r]^{\phantom{abcd}i'_*}\ar@{=}[d]&H^{p}_{q,\nr}(X\times \mathbb{P}^1,n)\ar@{->>}[r]^{j'^*}\ar@{->>}[d]^{f^*}&H^{p}_{q,\nr}(X\times \mathbb{A}^1,n)\ar[r]^{\delta=0\phantom{abcd}}\ar[d]^{f^*}&H^{p-1}_{q,\nr}(X\times \{0\},n-1)\ar@{=}[d]\\
H^{p-2}_{q-1,\nr}(X\times \{0\},n-1)\ar[r]^{\phantom{abcd}i_*}&H^{p}_{q,\nr}(X\times \mathbb{A}^1,n)\ar[r]^{j^*}&H^{p}_{q,\nr}(X\times \mathbb{G}_m,n)\ar[r]^{\delta\phantom{abcd}}&H^{p-1}_{q,\nr}(X\times \{0\},n-1),
}
\end{equation*}  
\end{small}

where $f\colon X\times (\mathbb{P}^1\setminus \{\infty\})\hookrightarrow X\times \mathbb{P}^1$ is an open immersion. It suffices to show $i_*$ is a trivial map by \Cref{prop-A1 homotopy invariance}. By \Cref{{computation trivial proj bundle}}, the exact sequence 
\begin{equation*}
   0\to H^{p-2}_{q-1,\nr}(X\times Q,A(n-1))\xrightarrow{g^*} H^{p}_{q,\nr}(X\times \mathbb{P}^1,A(n))\xrightarrow{h^*} H^{p}_{q,\nr}(X,A(n))\to 0
\end{equation*}
is split, where $g^*$ is induced by closed embedding of a rational point $Q\hookrightarrow \mathbb{P}^1$, $h^*$ is the open restriction, and $h^*\circ \pi'^*=\id$ for $\pi'\colon X\times \mathbb{P}^1\to X$. So we can find a canonical projection (i.e., independent of the choice of the rational point $Q$) $$p\colon H^{p}_{q,\nr}(X\times \mathbb{P}^1,A(n))\to H^{p-2}_{q-1,\nr}(X\times Q,A(n-1))= H^{p-2}_{q-1,\nr}(X,A(n-1)) $$ such that $p\circ g^*=\id$. Then we have 
\begin{equation*}
\ker(i_*)=\ker(f^*\circ i'_*)=p(\ker (f^*))\subset H^{p-2}_{q-1,\nr}(X\times \{0\},n-1). 
\end{equation*}
Now, take $Q=\infty$ and take $p\colon H^{p}_{q,\nr}(X\times \mathbb{P}^1,A(n))\to H^{p-2}_{q-1,\nr}(X\times \{\infty\},A(n-1))$. Then $p(\ker(f^*))=\im(p\circ g^*)=H^{p-2}_{q-1,\nr}(X\times \{\infty\},A(n-1))=H^{p-2}_{q-1,\nr}(X\times \{0\},A(n-1))$, as we whised.
\end{proof}

\section{Long Exact Sequence Involving Higher Chow Groups}
In this section, we prove that Bloch's higher cycle class map (see \cite[Section 4]{bloch11986algebraic}) with finite coefficients of a quasi-projective equi-dimensional scheme over a field fits naturally in a long exact sequence involving refined unramified cohomology. We restrict ourselves to the Borel--Moore cohomology theories (a) and (d) from \Cref{section-2.1} for convenience. We write $H^i_{\e t}(X,M(n))$ for the \'etale Borel--Moore cohomology with finite coefficients $M(n)\coloneqq\Z/m(n)=\mu_m^{\otimes n}$; we also write $H_i^{BM}(X,M(n))\coloneqq H^{2\dim(X)-i}(X,M(\dim(X)-n))$ for the \'etale Borel--Moore homology. Moreover, let $Z_i(X,\ast)$ be Bloch's cycle complex defined in \cite{bloch1986algebraic}, and we write $\CH_i(X,n;A)$ for Bloch's higher Chow groups with coefficients an abelian group $A$, that is 
\[
\CH_i(X,n;A)\coloneqq H_n(Z_i(X,\ast)\otimes^L A),
\]
where for the case $A=\Z$ we write $\CH_i(X,n)\coloneqq \CH_i(X,n;\Z)$. We recommend the article \cite{bloch1986algebraic,bloch11986algebraic} to readers regarding the definition and basic properties of higher Chow groups. The higher Chow groups are a motivic homology theory, so we have a spectral sequence associated to the `coniveau' filtration (indeed it is a niveau filtration).

\begin{thm}[{\cite[Section 10]{bloch1986algebraic}}]\label{thm-spectral sequence for higher chow}
 For any equi-dimensional algebraic $k$-scheme $X$, there is a spectral sequence 
 \begin{equation}\label{equ-ss for higher chow}
  E_1^{p,q}(b)\coloneqq \bigoplus_{x\in X^{(p)}}\CH^{b-p}(x,-p-q)\Rightarrow N^{\bullet}\CH^b(X,-p-q).
 \end{equation}
Here $\CH^{b-p}(x,-p-q)\coloneqq \varinjlim\limits_{U\subset \overline{\{x\}}}\CH^{b-p}(U,-p-q)$ where $U$ runs through all non-empty open subsets of the closure of $x\in X$, denoted by $\overline{\{x\}}$, and $$N^h\CH^b(X,n):=\varinjlim\limits_{\codim (Z)=c\geq h}\im(\CH^{b-c}(Z,n)\to \CH^b(X,n))$$ where $Z$ runs through all closed subsets of $X$ with $\codim (Z)\coloneqq\dim(X)-\dim(Z)=c\geq h$.
\end{thm}

 For smooth equi-dimensional $X$ the higher cycle class map
 $$\CH^b(X,n)\to \CH^b(X,n;M)\xrightarrow{\cl_M^{b,n}} H^{2b-n}_{\e t}(X,M(b))$$
  defined by Bloch (see\cite[Section 4]{bloch11986algebraic}) is compatible with localization sequences for any smooth closed subset of $X$(see \Cref{prop-higher cycle class map for any field}), which implies that  the higher cycle class map $\CH^b(X,n)\to H^{2b-n}_{\e t}(X,M(b))$ is compatible with the coniveau filtrations starting from $E_1$ for all smooth schemes.
 
Moreover, by \cite[Theorem 10.1]{bloch1986algebraic} for smooth equi-dimensional $X$, the groups $\CH^b(X,-q)$ satisfy Gersten's conjecture: the complex of Zariski sheaves (cf. \cite[(10.0.2)]{bloch1986algebraic})
 \begin{equation*}
     \bigoplus_{x\in X^{(0)}}(i_x)_\ast\CH^b(x,-q)\to \bigoplus_{x\in X^{(1)}}(i_x)_\ast\CH^{b-1}(x,-q-1)\to \dots \to \bigoplus_{x\in X^{(-q)}}(i_x)_\ast\CH^{b+q}(x,0)\to 0
 \end{equation*}
  is a flasque resolution of the Zariski sheaf $\mathcal{CH}^b_X(-q)$; the sheaf associated to the Zariski presheaf $U\mapsto \CH^b(U,-q)$. Furthermore, the $E_2^{p,q}(b)$ is precisely the second page of the local to global spectral sequence computing higher Chow groups; see \cite[Section 10]{bloch1986algebraic}.

  
For higher Chow groups with finite coefficients $M$, we also have the following consequence of the work of Geisser--Levine \cite{geisser2001bloch} and the Bloch–-Kato conjecture, which is proven by Voevodsky \cite[Theorem 6.1]{voevodsky2011motivic}.
  
\begin{prop}[cf.~{\cite[Theorem 1.1]{geisser2001bloch},\cite[Theorem 9.1]{kerz2012cohomological}}]\label{cor-GL01}
 Let $m$ be a fixed integer invertible in a field $k$, and let $X$ be an essentially smooth quasi-projective $k$-scheme (e.g. $X$ is a localization of a smooth scheme of finite type over $k$). Bloch's cycle class map $\CH^b(X,n;\Z/m\Z)\to H_{\e t}^{2b-n}(X,\mu_m^{\otimes b})$ is an isomorphsim for $b\leq n$ and is an injection for $b=n+1$. In particular, for every point $x\in X$, the cycle class map $\CH^b(x,n;\Z/m\Z)\to H_{\e t}^{2b-n}(x,\mu_m^{\otimes b})$ is an isomorphsim for $b\leq n$.
\end{prop}
\begin{proof}
 By referring to \cite[Theorem 1.1]{geisser2001bloch} and based on Bloch--Kato conjecture \cite{voevodsky2011motivic}, the following natural map in the derived category $D(X_{\rm Zar},\Z)$ is a quasi-isomorphism for all $b$,
 $${\rm cl}_X^b\colon \Z/m\Z(b) \to \tau^{\leq b}R\pi_*(\mu_m^{\otimes b}),$$ 
 where
 \begin{itemize}
     \item $\Z/m\Z(b)$ is the motivic complex of weight $b$ with coefficient $\Z/m\Z$, given by 
     \[
     \Z/m\Z(b)\coloneqq \Z(b)\otimes^L \Z/m\Z\coloneqq C_*\Z_{tr}(\mathbb{G}_m^{\wedge b})[-b]\otimes^L \Z/m\Z \text{\quad(cf.~\cite{mazza2006lecture})};
     \]
     \item $\pi\colon X_{\mathrm {\e t}}\to X_{\rm Zar}$ is the identity continuous morphism of sites;  
     \item $\tau^{\leq b}R\pi_*(\mu_m^{\otimes b})$ is the $b$-th canonical truncation of the derived complex $R\pi_*(\mu_m^{\otimes b})$. 
 \end{itemize}
Bloch's higher cycle class map $\CH^b(X,n;\Z/m\Z)\to H_{\e t}^{2b-n}(X,\mu_m^{\otimes b})$ is induced by the following morphism in the derived category $D(X_{\rm Zar},\Z)$
\begin{align}\label{new morphism}
{\rm cl}_X^b\colon \Z/m\Z(b) \xrightarrow{\sim} \tau^{\leq b}R\pi_*(\mu_m^{\otimes b})\to R\pi_*(\mu_m^{\otimes b}),
\end{align}
by \cite{Voevodsky2002motivic} and \cite{geisser2001bloch}, where the last arrow is the canonical truncated morphism. Using the following spectral sequence
\[
E^{p,q}_2=H^p_{\rm Zar}(X,R^q\tau^{\geq b+1}R\pi_*(\mu_m^{\otimes b}))\Rightarrow \mathbb{H}^{p+q}_{\rm Zar}(X,\tau^{\geq b+1}R\pi_*(\mu_m^{\otimes b})),
\]
where
\begin{equation*}
 E_2^{p,q}= 
	\begin{cases}
	0 & \text{if $q<b+1$};\\
	  H^p_{\rm Zar}(X,R^q\pi_*(\mu_m^{\otimes b})) & \text{if $q\geq b+1$},
	\end{cases}   
\end{equation*}
 we conclude that $\mathbb{H}^{p+q}_{\rm Zar}(X,\tau^{\geq b+1}R\pi_*(\mu_m^{\otimes b}))=0$ if $p+q\leq b$. So taking the Zariski hyper-cohomology of the morphism (\ref{new morphism}), we conclude that $H^{n}_M(X,\Z/m\Z(b))\to \mathbb{H}^{n}(X,\tau^{\leq b}R\pi_*(\mu_m^{\otimes b}))$ is an isomorphism for $n\leq b$ and is an injection for $n=b+1$. Equivalently, Bloch's cycle class map $\CH^b(X,n;\Z/m\Z)\to H_{\e t}^{2b-n}(X,\mu_m^{\otimes b})$ is an isomorphsim for $b\leq n$ and an injection for $b=n+1$. Finally, for every $x\in X$ and for all $b\leq n$, we conclude our proof from
 \begin{equation*}
  \CH^{b}(x,n;\Z/m\Z)\coloneqq \varinjlim\limits_{U\subset \overline{\{x\}}}\CH^{b}(U,n;\Z/m\Z)\xrightarrow{\sim }\varinjlim\limits_{U\subset \overline{\{x\}}} H^{2b-n}(U,\mu_m^{\otimes b})\eqqcolon H_{\e t}^{2b-n}(x,\mu_m^{\otimes b}).\qedhere
 \end{equation*} 
\end{proof}
\subsection{Higher Cycle Class Map}\label{construction higher cycle class map}
Let $m$ be an integer in $k^*$ and $M\coloneqq \Z/m\Z$.
In this section we give a construction of the higher cycle class map $$\cl^M_{i,n}\colon\CH_i(X,n;M)\to H^{BM}_{2i+n}(X,M(i))$$ for every quasi-projective equi-dimensional $k$-scheme $X$ and show it is compatible with the localization long exact sequences. This is probably well-known to experts, but the authors could not find a precise reference. We write $f\colon X\to \mathrm{Spec}(k)$ for the structure morphism in this section.

\begin{prop}[{cf. \cite[Corollary 8.12]{cisinski2014integral},\cite[(11)]{kelly2014isomorphisme}}]\label{higher Chow motivic isom}
Let $X$ be a quasi-projective equi-dimensional scheme over a perfect field $k$.
Then for all $i,j\in \Z$, there are canonical isomorphisms:
\begin{equation*}
  \CH_i(X,j-2i;M)\cong \left\{\begin{array}{ll}
     Hom_{DM_{cdh}(k,M)}(M(i)[j],Rf_\ast f^!M) & \text{if } i\geq 0 \\
      Hom_{DM_{cdh}(k,M)}(M,Rf_\ast f^!M(-i)[-j]) & \text{if } i\leq0,
  \end{array} \right.
\end{equation*}
where $DM_{cdh}(k,M)$ is the category of  $cdh$-motives over Spec$(k)$ and coefficients in $M$(cf.~\cite[Lecture 9]{mazza2006lecture}, \cite{cisinski2014integral}), and $M(i)\coloneqq \Z(i)\otimes^LM$.
\end{prop}
\begin{proof}
The proof is essentially the same as \cite[Chapter 5 Proposition 4.2.9]{voevodsky2000cycles}. Let us consider first the case $i\geq 0$. We have the following canonical isomorphisms

\begin{equation*}
    \begin{split}
     Hom_{DM_{cdh}(k,M)}(M(i)[j],Rf_\ast f^!M)&\cong A_{i,j-2i}({\rm Spec}(k),X)_M \quad \text{(by \cite[Theorem 8.11]{cisinski2014integral})}\\
     \text{(by \cite[(8.3.1)]{cisinski2014integral} }\quad&=\mathbb H_{cdh}^{-(j-2i)}(\mathrm{Spec}(k),\underline C_\ast(z_{equi}(X,i))_{cdh}\otimes^LM)\\
   \text{(by \cite[Theorem 5.13]{SV00})}\quad &\cong H_{j-2i}((\underline C_\ast(z_{equi}(X,i))_{cdh}\otimes^LM)({\rm Spec}(k)))\\
    \text{(by definition)}\quad&=H_{j-2i}(C_\ast(z_{equi}(X,i))\otimes^LM)\\
    &\cong \CH_i(X,j-2i;M).
    \end{split}
\end{equation*}

Here the last isomorphism is induced by the quasi-isomorphsim of complexes $$C_\ast(z_{equi}(X,i))\to Z_{i}(X,*)$$ if $k$ is of characteristic zero (see \cite[Chapter 5 Proposition 4.2.9]{voevodsky2000cycles}), and is induced by the quasi-isomorphsim of complexes $$C_\ast(z_{equi}(X,i))[\frac{1}{p}]\to Z_{i}(X,*)[\frac{1}{p}]$$ if $k$ is of characteristic $p>0$ (see \cite[Théorème 2.3]{kelly2014isomorphisme}).

Suppose now that $i\leq 0$. In this case, we have the following canonical isomorphisms
\begin{equation*}
    \begin{split}
        \CH_i(X,j-2i;M)&\cong \CH_0(X\times \mathbb{A}^{-i},j-2i;M) \quad \text{(by \cite[Theorem 2.1]{bloch1986algebraic})}\\
     \text{(by case $i\geq 0$)}\quad   &\cong Hom_{DM_{cdh}(k,M)}(M(0)[j-2i],Rg_\ast g^!M)\\
        &\cong Hom_{DM_{cdh}(k,M)}(M(0)[j-2i],Rf_\ast f^!M(-i)[-2i])\\
        \quad &\cong Hom_{DM_{cdh}(k,M)}(M,Rf_\ast f^!M(-i)[-j]),
    \end{split}
\end{equation*}
where $g:X\times \mathbb{A}^{-i}\to {\rm Spec}(k)$ is the structure morphism. Moreover, the third isomorphism is given by the isomorphisms 
\[
Rg_\ast g^!M\cong M^c(X\times \mathbb{A}^{-i})\otimes^LM\cong M^c(X)(-i)[-2i]\otimes^LM\cong Rf_\ast f^!M,
\]
where the first (and the last) isomorphism is induced by \cite[Proposition 8.10]{cisinski2014integral} and $M^c(*)$ is the \emph{motive with compact support} (cf.~\cite[Lecture 16]{mazza2006lecture}), and the middle one is induced by homotopy invariance (see \cite[Chapter 5 Corollary 4.18]{voevodsky2000cycles} if $k$ is of characteristic zero and see \cite[Corollary 5.5.9]{kellyphdthesis} if $k$ is of characteristic $p>0$).
\end{proof}

In \cite[Corollary 9.5]{cisinski2014integral} it is shown that there is a realization functor for every noetherian scheme $X$
\[
(-)_h\colon DM_{cdh}(X,M)\to DM_h(X,M),
\]
which preserves the six operations when restricted to constructible objects by \cite[Remark 9.6]{cisinski2014integral}. Using this together with the fact that there is a canonical equivalence for finite $M$ (\cite[Remark 9.6]{cisinski2014integral}, \cite[Corollary 5.5.4]{CisinskiDenis-Charles2016Em}) 
\begin{equation}\label{equ-an equivalence}
    DM_h(X,M)\cong D(X_{\e t},M)
\end{equation}
where $D(X_{\e t},M)$ is the derived category of the abelian category of sheaves of $M$-modules on the small \'etale site of $X$, we obtain the cycle class map as follows.
\begin{define}\label{def-realization}
Let $X$ be a quasi-projective equi-dimensional scheme over a perfect filed $k$. 
The higher cycle class map $\cl_{i,n}^M\colon\CH_i(X,n,M)\to H^{BM}_{2i+n}(X,M(i))$ for $i\geq 0$ is defined as the following composition
\begin{equation}\label{equ-construction of higher cycle class mapI}
 \begin{split}
     \CH_i(X,n,M)&\cong Hom_{DM_{cdh}(k,M)}(M(i)[2i+n],Rf_\ast f^!M) \quad \text{(by \Cref{higher Chow motivic isom})}\\
     &\xrightarrow{\text{h-sheafification}} Hom_{DM_h(k,M)}(M(i)[2i+n],Rf_\ast f^!M)\\
     &\cong Hom_{D_{\e t}(Spec(k),M)}(M(i)[2i+n],Rf_\ast f^!M)\quad\text{by (\ref{equ-an equivalence})}\\
     &\eqqcolon H^{BM}_{2i+n}(X,M(i)).
 \end{split}   
\end{equation}
Similarly, for $i\leq 0$, we define the higher cycle class map $\cl_{i,n}^{M}$ as the following composition
\begin{equation}\label{equ-construction of higher cycle class mapII}
 \begin{split}
     \CH_i(X,n,M)&\cong Hom_{DM_{cdh}(k,M)}(M,Rf_\ast f^!M(-i)[-2i-n]) \quad \text{(by \Cref{higher Chow motivic isom})}\\
     &\xrightarrow{\text{h-sheafification}} Hom_{DM_h(k,M)}(M,Rf_\ast f^!M(-i)[-2i-n])\\
     &\cong Hom_{D_{\e t}(Spec(k),M)}(M,Rf_\ast f^!M(-i)[-2i-n])\quad\text{by (\ref{equ-an equivalence})}\\
     &\cong Hom_{D_{\e t}(Spec(k),M)}(M(i)[2i+n],Rf_\ast f^!M)\\
     &\eqqcolon H^{BM}_{2i+n}(X,M(i)).
 \end{split}   
\end{equation}

For equi-dimensional smooth $X$ of dimension $d$, we also obtain the cycle class map $$\cl_M^{i,n}\colon\CH^i(X,n;M)=\CH_{d-i}(X,n;M)\xrightarrow{\cl_{d-i,n}^M} H^{BM}_{2d-2i+n}(X,M(d-i))\cong H_{\e t}^{2i-n}(X,M(i))$$
by Poincar\'e duality.
\end{define}

\begin{prop}\label{prop-cycle class map is compatible with localization}
The higher cycle class map defined in \Cref{def-realization} is compatible with flat pullback and proper pushforward. Moreover, it respects the localization long exact sequence.
\end{prop}
\begin{proof}
The first statement follows from the fact that the realization functor preserves the six operations (see \cite[Corollary 9.5]{cisinski2014integral}), which determine the morphisms. 

For the second statement we also use the compatibility with the six operations, together with \cite[Theorem 5.11]{cisinski2014integral} which says that $DM_{cdh}(-,R)$ has the localization property for any $\Z[\frac{1}{p}]$-algebra $R$, and the proof will then again be similar to \cite[Chapter 5 Proposition 4.2.9]{voevodsky2000cycles}. This namely implies that if $Z\subseteq X$ is a closed subscheme with open complement $U$, we get a distinguished triangle 
\[
R(f|_Z)_\ast (f|_Z)^!\dagger\to Rf_\ast f^!\dagger\to R(f|_U)_\ast (f|_U)^!\dagger\overset{+1}{\to}
\]
for every $\dagger\in DM_{cdh}(k,M)$. The realization functor maps this to the localization sequence in $D(M)$. As the morphisms involved in the construction of the cycle class map are all functorial, this gives the result.
\end{proof}

\begin{prop}\label{BL conjecture}
Let $X$ be a smooth quasi-projective scheme over a perfect field $k$.
 The cycle class map $\cl^M_{i,j-2i}\colon \CH_i(X,j-2i,M)\to H_{\e t}^{2d-j}(X,M(d-i))$ is (up to natural isomorphism) the same as the one from \cite{geisser2001bloch}. In particular for such $X$ the map $\cl_M^{i,n}$ is an isomorphism if $n\geq i$ and injective for $n=i-1$.
\end{prop}
\begin{proof}
We cover the case where $i\geq 0$, the $i\leq 0$ case is similar. By \cite[Corollary 5.9]{cisinski2014integral} and \Cref{higher Chow motivic isom} we know that 
\begin{equation*}
    \begin{split}
        \CH_i(X,j-2i;M)&\cong Hom_{DM_{cdh}(k,M)}(M(i)[j],M^c(X))\\
        &\cong Hom_{DM(k,M)}(M(i)[j],M^c(X))\\
        &=Hom_{DM^{eff}(k,M)}(M(i)[j],M^c(X)),
    \end{split}
\end{equation*}
where $DM(k,M)$ is the category of motives over $k$ with $M$-coefficients and $DM^{eff}(k,M)$ is the category of effective  motives over $k$ with $M$-coefficients (cf. \cite[Lecture 14]{mazza2006lecture}). And the last equation holds because of the fully faithfulness of $DM^{eff}(k,M)\to DM(k,M)$ by the cancellation theorem \cite[Corollary 8.2]{cisinski2014integral}. Thus, $\alpha^*:DM^{eff}(k,M)\to DM^{eff}_{\e t}(k,M)$ (i.e., the \'etale-sheafification functor) induces a morphism 
\begin{equation}\label{equ-higher cycle class map for smooth schemes}
   \CH_i(X,j-2i;M)\xrightarrow{\alpha^*} Hom_{DM^{eff}_{\e t}(k,M)}(M(i)[j],\alpha^\ast M^c(X))\cong H^{2d-j}_{\e t}(X,M(d-i); 
\end{equation}
here $DM^{eff}_{\e t}(k,M)$ is category of effective \'etale motives over $k$ with $M$-coefficients (see \cite[Lecture 9]{mazza2006lecture}) and the isomorphism is given by \cite[Proposition 4.3; Lemma~4.4]{kelly2017weight}.


By \cite[Section 3]{ivorra2010realisation} (cf. \cite[Section 4]{kelly2017weight}), the morphism (\ref{equ-higher cycle class map for smooth schemes}) is exactly the higher cycle class map defined in \cite{geisser2001bloch}. Now we check that this is the same map as we defined in \Cref{def-realization}. By \cite[Theorem 4.5.2, Theorem 5.5.3]{CisinskiDenis-Charles2016Em}, the isomorphism $D(k_{\e t},M)\cong DM_h(k,M)$ used in the construction of the higher cycle class map before is defined as the composition of the canonical functors
\[
D(k_{\e t},M)\xrightarrow[\rho_!]{\sim}DM_{\e t}^{eff}(k,M)\overset{\sim}{\to}DM_h^{eff}(k,M)\overset{\sim}{\to}DM_h(k,M),
\]
where $\rho_!\colon D_{\e t}(M)\overset{\sim}{\to} DM_{\e t}^{eff}(k,M)$ is from \cite[Theorem 4.5.2]{CisinskiDenis-Charles2016Em}. We summarize the maps involved in the following diagram

\[
\begin{tikzcd}
DM^{eff}(k,M)\ar[r,"\text{f.f.}"]\ar[d,"\alpha^\ast"]& DM(k,M)\ar[d]\ar[ddr,"\cong"]&\\
DM^{eff}_{\e t}(k,M)\ar[r,"\cong"]\ar[d,"\cong"]& DM_{\e t}(k,M)\ar[d,"\cong"]&\\
DM^{eff}_h(k,M)\ar[r,"\cong"]& DM_h(k,M)& DM_{cdh}(k,M)\nospacepunct{.}\ar[l,swap,"h"]
\end{tikzcd}
\]
As all the maps are induced by the sheafification functors on sites, one can check that this diagram commutes (up to isomorphism) and preserve the six operations. Now, on the one hand, the higher cycle class maps in \Cref{def-realization} is defined as the compositions
\begin{equation*}
    DM_{cdh}(k,M)\xrightarrow{h}DM_h(k,M)\xleftarrow{\sim}DM_h^{eff}(k,M)\xleftarrow{\sim}DM_{\e t}^{eff}(k,M)\xleftarrow[\rho_!]{\sim} D(k_{\e t},M);
\end{equation*}
on the other hand, the map (\ref{equ-higher cycle class map for smooth schemes}) is defined as 
\begin{equation*}
  DM_{cdh}(k,M)\xleftarrow{\sim}  DM(k,M)\xleftarrow{\text{f.f.}}DM^{eff}(k,M)\xrightarrow{\alpha^*}DM_{\e t}^{eff}(k,M)\xleftarrow[\rho_!]{\sim} D(k_{\e t},M).
\end{equation*}
The above commutative diagram implies that both cycle class maps constructed above coincide up to natural isomorphism. The second statement is now \Cref{cor-GL01}.
\end{proof}

\begin{cor}\label{cor-cycle class map is an iso}
Let $X$ be a quasi-projective equi-dimensional scheme of dimension $d$ over a perfect filed $k$.  The cycle class map $\cl^M_{i,n}\colon \CH_i(X,n;M)\to H^{BM}_{2i+n}(X,M(i))$ is an isomorphism for all $i\geq d-n$. Moreover, $\cl^M_{d-n-1,n}$ is an injection for such~$X$.
\end{cor}
\begin{proof}
We do induction by the dimension of $X$. If $d=0$, then it is automatically smooth and hence it follows from \Cref{BL conjecture} or \Cref{cor-GL01}. Moreover, we may assume that $X$ is reduced as $\CH_i(X,n;M)\cong \CH_i(X_{\mathrm{red}},n;M)$ and $H^{BM}_{2i+n}(X,M(i))\cong H^{BM}_{2i+n}(X_{\mathrm{red}},M(i))$ by localization sequences. Let $Y\subset X$ be an equi-dimensional closed subscheme with $X^{\mathrm{sing}}\subset Y\subset X$, where $X^{\mathrm{sing}}$ is the singular locus of $X$; let $\emptyset \neq U\subset X$ be the open complement of $Y$. Then $U$ is smooth, and for all $i\geq d-n$ we have a commutative diagram between localization sequences by \Cref{prop-cycle class map is compatible with localization}
\begin{equation*}
  \adjustbox{scale= .85}{
\xymatrix{\CH_i(U,n+1;M)\ar[d]\ar[r]&\CH_i(Y,n;M)\ar[d]\ar[r]&\CH_i(X,n;M)\ar[d]\ar[r]&\CH_i(U,n;M)\ar[d]\ar[r]&\CH_i(Y,n-1;M)\ar[d]\\
H^{BM}_{2i+n+1}(U,M(i))\ar[r]&H^{BM}_{2i+n}(Y,M(i))\ar[r]&H^{BM}_{2i+n}(X,M(i))\ar[r]&H^{BM}_{2i+n}(U,M(i))\ar[r]&H^{BM}_{2i+n-1}(Y,M(i))
}  }  
\end{equation*} 
where all the outer vertical maps are higher cycle class maps, so are isomorphisms by induction and \Cref{BL conjecture}. Hence the five-lemma implies that the middle map is an isomorphism as well. Similarly, $\cl^M_{d-n-1,n}$ is an injection for every quasi-projective equi-dimensional $X$.
\end{proof}
Now we extend our construction to non-perfect fields.
\begin{lem}[{cf.~\cite[Lemma 6.8]{schreieder2021refined}}]\label{lem-purely inseparable extension}
 Let $X$ be a separated scheme of finite type over a field $k$ of characterstic $p> 0$ and let $E/k$ be a purely inseparable extension. Let $f:X_E\to X$ be the canonical morphism. Then the flat pullback $f^*:\CH_i(X,n;M)\to \CH_i(X_E,n;M)$ is a canonical isomorphism.
\end{lem}
\begin{proof}
It suffices to consider the case where $E/k$ is a finite extension of degree $p^s$ for some $s$ as $\CH_i(X_E,n;M)=\varinjlim\limits_{F}\CH_i(X_F,n;M)$ and $F$ ranges over all finitely generated field extensions of $k$ contained in $E$. Then we have $f_*\circ f^*=p^s\cdot\id$ by definition (see \cite[Corollary 1.4]{bloch1986algebraic}). Note that $p^s\cdot\id:\CH_i(X,n;M)\xrightarrow{\sim}\CH_i(X,n;M)$ is an isomorphism as $p$ is invertible in $M=\Z/m\Z$,  so $f^*$ is injective. Moreover, by the same argument as in \cite[Lemma 6.8]{schreieder2021refined}, we have $f^*\circ f_*=p^t\cdot\id$ for some $t$ and so $f^*$ is surjective.
\end{proof}
\begin{prop}\label{prop-higher cycle class map for any field}
Let $X$ be a quasi-projective equi-dimensional scheme of dimension $d$ over a field $k$. Then there exists a canonical higher cycle class map for all $i,n\in \Z$ $$\cl_{i,n}^M\colon\CH_i(X,n,M)\to H^{BM}_{2i+n}(X,M(i))$$
satisfying the following:
\begin{itemize}
    \item[(i)] it is compatible with localization exact sequences, flat pullback and proper push-forward;
    \item[(ii)] it is an isomorphism for all $i\geq d-n$ and is an injection for $i=d-n-1$;
\item[(iii)] for smooth $X$, it is the same as the one from \cite{geisser2001bloch}.
\end{itemize}
\end{prop}
\begin{proof}
 For perfect fields, we have done this in the beginning of this subsection (i.e., \Cref{def-realization}, \Cref{prop-cycle class map is compatible with localization}, \Cref{BL conjecture} and \Cref{cor-cycle class map is an iso}). So let us assume that $k$ is non-perfect, and let $k'$ be the perfect closure of $k$. Then $k'/k$ is a purely inseparable extension. Let $f:X_{k'}\to X$ be the canonical morphism, and then we define the higher cycle class map for $X$ as the following
 \begin{small}
 \begin{equation}
     \CH_i(X,n;M)\xrightarrow[\sim]{(f^*)^{-1}}\CH_i(X_{k'},n;M)\xrightarrow{\cl_{i,n}^M(X_{k'})}H_{2i+n}^{BM}(X_{k'},M(i))\xrightarrow[\sim]{(f^*)^{-1}}H_{2i+n}^{BM}(X,M(i)),
 \end{equation}
 \end{small}
 
 where the first and the last arrows are the natural isomorphisms from \Cref{lem-purely inseparable extension} and the topological invariance of \'etale cohomology. Therefore, $\cl_{i,n}^M(X)$ satisfies (i), (ii) and (iii) just like $\cl_{i,n}^M(X_{k'})$ does.
\end{proof}
\subsection{Refined Higher Chow Groups}
Here we apply a construction to the higher Chow groups, analogous to the construction of the refined unramified cohomology groups. For this we put $$\CH^i(F_pX,n)\coloneqq \varinjlim\limits_{F_pX\subseteq U\subseteq X}\CH^i(U,n).$$ As the higher Chow groups satisfy the localization sequence \cite[Theorem 3.1]{bloch1986algebraic}, we get analogous to \Cref{prop-Gysin sequence} the following.

\begin{prop}\label{prop-Gysin sequence higher Chow}
For any $b,n,p,m\geq 0$, there is a long exact sequence for every algebraic scheme $X$
\[
\xrightarrow{\iota_\ast}\CH^b(F_{p+m}X,n)\to\CH^b(F_{p-1}X,n)\xrightarrow{\partial}\varinjlim_{\codim(Z)=p}\CH^{b-p}(F_mZ,n-1)\xrightarrow{\iota_\ast}\CH^b(F_{p+m}X,n-1)\to,
\]
where $Z$ runs through all closed subset of $X$ satisfying $\codim(Z)\coloneqq\dim(X)-\dim(Z)=p$.
\end{prop}
\begin{rem}
 \rm{Applying \Cref{thm-spectral sequence for higher chow} and  \Cref{prop-Gysin sequence higher Chow}, we have 
 \begin{equation*}
    \begin{split}
      N^h\CH^b(X,n)&:=\varinjlim\limits_{\codim (Z)=c\geq h}\im(\CH^{b-c}(Z,n)\to \CH^b(X,n))\\
      &\cong  \ker(\CH^b(X,n)\to \CH^b(F_{h-1}X,n)),
    \end{split} 
 \end{equation*}
 which is the same filtration that used in \cite[(5.1)]{schreieder2021refined}. However, for $n=0$ this filtration is different from the filtration in \cite[Definition 7.3]{schreieder2021refined} as here $N^0\CH^b(X)=\CH^b(X)$ while $N^0\CH^b(X)=\CH^b(X)_{\rm hom}$ for the filtration in loc.~cit.
}
\end{rem}
Taking $m=0$ in \Cref{prop-Gysin sequence higher Chow} gives the following long exact sequence.
\begin{prop}\label{prop-Gysin sequence higher Chow generic point}
There is a long exact sequence for every algebraic $k$-scheme $X$
\[
\to \bigoplus_{x\in X^{(p)}}\CH^{b-p}(\kappa(x),n)\xrightarrow{\iota_\ast} \CH^b(F_pX,n)\to\CH^b(F_{p-1}X,n)\xrightarrow{\partial}\bigoplus_{x\in X^{(p)}}\CH^{b-p}(\kappa(x),n-1)\xrightarrow{\iota_\ast}.
\]
\end{prop}

As a direct consequence we obtain the analogous statement of \Cref{cor-degeneration of refined unramified cohomology}.
\begin{cor}\label{cor-degeneration of higher Chow groups}
Let $X$ be an algebraic $k$-scheme. We have natural isomorphisms $\CH^b(X,n)\xrightarrow{\sim}\CH^b(F_pX,n)$ for $p\geq b$ and if $0\leq p < b-n$ we have $\CH^b(F_pX,n)\xrightarrow{\sim}\CH^b(F_0X,n)=0$.
\end{cor}

\subsection{Long Exact Sequence}\label{section-4.3 long exact sequence}
In this section, we are going to show the relation between the (co)kernel of the higher cycle class map and refined unramified cohomology. Before we begin the proof of \Cref{thm-main thm2} (=\Cref{thm-our main intro}), let us first prove the following lemma.

\begin{lem}\label{lem-a claim}
Let $k$ be a field and let $b$ and $n$ be two integers with $b>n\geq 0$. Suppose $X$ is a scheme over a field $k$ satisfying either of the following two conditions: 
\begin{itemize}
    \item[(a)] $X$ is a quasi-projective equi-dimensional $k$-scheme;
    \item[(b)] there exists an open subset $U\subset X$ with $\dim(X\setminus U)\leq \dim(X)-(b-n+1)$ such that $U$ is a smooth $k$-scheme.
\end{itemize}
Let $m$ be an integer invertible in $k$. Then the following is exact 
\begin{small}
\begin{equation}\label{equ-a claim}
 \CH^b(F_{b-n}X,n;\Z/m\Z)\xrightarrow{\cl} H^{2b-n}(F_{b-n}X,\Z/m\Z(b))\xrightarrow{\rm restr.}H^{2b-n}(F_{b-n-1}X,\Z/m\Z(b)),
\end{equation}
\end{small}

where the morphism `$\cl$' is given by \Cref{prop-higher cycle class map for any field} if $X$ satisfies (a) and is given by \cite[Section 4]{bloch11986algebraic} if $X$ satisfies (b).
\end{lem}
\begin{proof}
Note that we have the following exact sequence by \Cref{prop-Gysin sequence higher Chow generic point}
\begin{equation*}
    \bigoplus_{x\in X^{(b-n)}}\CH^n(x,n;\Z/m\Z)\to \CH^b(F_{b-n}X,n;\Z/m\Z)\to\CH^b(F_{b-n-1}X,n;\Z/m\Z)=0,
\end{equation*}
where $\CH^b(F_{b-n-1}X,n;\Z/m\Z)\cong \CH^b(F_0X,n;\Z/m\Z)=0$ by \Cref{cor-degeneration of higher Chow groups} as $b>n$ and $b-n-1<b-n$. Thus $\bigoplus\limits_{x\in X^{(b-n)}}\CH^n(x,n;\Z/m\Z)\to \CH^b(F_{b-n}X,n;\Z/m\Z)$ is surjective. On the other hand, by \cite[Theorem 9.1]{kerz2012cohomological} or \Cref{cor-GL01}, for every point $x\in U\subset X$, the cycle class map $\CH^n(x,n;\Z/m\Z)\cong H^n(x,\Z/m\Z(n))$ is an isomorphism. 
So we have the following equations
\begin{equation*}
\begin{split}
&\ker(H^{2b-n}(F_{b-n}X,\Z/m\Z(b))\xrightarrow{\rm restr.} H^{2b-n}(F_{b-n-1}X,\Z/m\Z(b)))\\
=&\im(\bigoplus_{x\in X^{(b-n)}}H^n(x,\Z/m\Z(n))\xrightarrow{\iota_*} H^{2b-n}(F_{b-n}X,\Z/m\Z(b)) )\\
=&\im(\CH^b(F_{b-n}X,n;\Z/m\Z)\xrightarrow{\cl} H^{2b-n}(F_{b-n}X,\Z/m\Z(b))),
\end{split}
\end{equation*}
where the first equation is given by \Cref{prop-Gysin sequence} and the last is induced by the following commutative diagram 
\begin{equation*}
\xymatrix{\bigoplus\limits_{x\in X^{(b-n)}}\CH^n(x,n;\Z/m\Z)\ar@{->>}[r]^{\iota_*}\ar[d]^{\sim}_{\cl}& \CH^b(F_{b-n}X,n;\Z/m\Z)\ar[d]_{\cl}\\
   \bigoplus\limits_{x\in X^{(b-n)}}H^n(x,\Z/m\Z(n))\ar[r]^{\iota_*}& H^{2b-n}(F_{b-n}X,\Z/m\Z(b))
    }
\end{equation*}
by applying \Cref{prop-cycle class map is compatible with localization} for case (a) and \cite[Section 4]{bloch11986algebraic} for case (b).
\end{proof}

\begin{prop}\label{prop-exact at H(X)}
 Let $X$ be a quasi-projective equi-dimensional scheme defined over a field $k$ and $m$ be an integer invertible in $k$.  The sequence
\[
\CH^b(X,n;\Z/m\Z)\xrightarrow{\rm cl} H^{2b-n}(X,\Z/m\Z(b))\xrightarrow{\rm restr.} H_{b-n-1,\nr}^{2b-n}(X,\Z/m\Z(b))
\]
is exact, where $\cl\colon\CH^b(X,n;\Z/m\Z)\to H^{2b-n}(X,\Z/m\Z(b))$ is the higher cycle class map from \Cref{construction higher cycle class map}.
\end{prop}
\begin{proof}
Taking finite coefficient $M\coloneqq\Z/m\Z$, the cycle class map $\mathrm{cl}^{b,n}$ is an isomorphsim if $b\leq n$ by \Cref{prop-higher cycle class map for any field}, so this sequence is exact as $H_{p,\nr}^{q}(X,\Z/m\Z(b))=0$ for any $p<0$. The sequence is also exact when $n=0$ by \Cref{lem-a claim} as $$H^{2b-n}(F_bX,\Z/m\Z(b))\cong H^{2b-n}(X,\Z/m\Z(b))$$ (see \Cref{cor-degeneration of refined unramified cohomology}), so we may assume that $b>n>0$. 

Consider the following commutative diagram 

\begin{small}
\begin{equation}\label{equ-diagram for main thm2}
    \xymatrix{    \varinjlim\limits_{\codim Z=b-n+1}\CH^{n-1}(Z,n;M)\ar[r]^{\sim}_{\cl}\ar[d]^{\iota_*}&\varinjlim\limits_{\codim Z=b-n+1}H^{n-2}(Z,M(n-1))\ar[d]^{\iota_*}&\\
    \CH^b(X,n;M)\ar[r]_{\cl}\ar[d]^{\rm restr.}&H^{2b-n}(X,M(b))\ar[r]^{\rm restr.}\ar[d]^{\rm restr.}&H^{2b-n}_{b-n-1,\nr}(X,M(b))\ar@{^(->}[d]\\
    \CH^b(F_{b-n}X,n;M)\ar[r]_{\cl}\ar[d]^{\partial}&H^{2b-n}(F_{b-n}X,M(b))\ar[r]^{\rm restr.}\ar[d]^{\partial}&H^{2b-n}(F_{b-n-1}X,M(b))\\
    \varinjlim\limits_{\codim Z=b-n+1}\CH^{n-1}(Z,n-1;M)\ar[r]^{\sim}_{\cl}\ar@{->>}[d]^{\iota_*}&\varinjlim\limits_{\codim Z=b-n+1}H^{n-1}(Z,M(n-1))\ar[d]^{\iota_*}\\
    \CH^b(X,n-1;M)\ar[r]_{\cl}&H^{2b-n+1}(X,M(b))\nospacepunct{,}&
    }
\end{equation}
\end{small}

where
\begin{itemize}
\item $Z$ runs through all the closed subsets of $X$ with $\dim(Z)=\dim(X)-(b-n+1)$;
\item the two  isomorphisms are the cycle class maps, based on \Cref{prop-higher cycle class map for any field};
\item the middle horizontal row is exact by \Cref{lem-a claim}.
\end{itemize}
Moreover, taking $j=b-n+1$ and $m\geq n-1$, then we have $$H^{2b-n}(F_{j+m}X,M(b))\cong H^{2b-n}(X,M(b))$$ and $$H^{2b-n+1}(F_{j+m}X,M(b))\cong H^{2b-n+1}(X,M(b))$$ by \Cref{cor-degeneration of refined unramified cohomology}, as we have  $j+m\geq b\geq \lceil(2b-n+1)/2\rceil\geq \lceil(2b-n)/2\rceil$. So applying \Cref{prop-Gysin sequence}, the second column in (\ref{equ-diagram for main thm2}) is exact. Similarly, applying \Cref{cor-degeneration of higher Chow groups} and \Cref{prop-Gysin sequence higher Chow}, the first column in (\ref{equ-diagram for main thm2}) is exact as well. Furthermore, the morphism $$\varinjlim\limits_{\codim Z=b-n+1}\CH^{n-1}(Z,n-1;M)\to \CH^b(X,n-1;M)$$ is surjective, since $\CH^b(F_{b-n}X,n-1;M)=0$ by \Cref{cor-degeneration of higher Chow groups}.


Now, we can prove \Cref{prop-exact at H(X)} directly by doing a diagram chase in (\ref{equ-diagram for main thm2}). We write down the details for the convenience of the reader. On the one hand, the composition 
\[
\CH^b(X,n;M)\to H^{2b-n}(X,M(b))\to H^{2b-n}_{b-n-1,\nr}(X,M(b))
\]
is zero as we have \Cref{lem-a claim} and $H^{2b-n}_{b-n-1,\nr}(X,M(b))\hookrightarrow H^{2b-n}(F_{b-n-1}X,M(b))$ is injective. Therefore, we have 
\[
\im(\CH^b(X,n;M)\to H^{2b-n}(X,M(b)))\subseteq \ker(H^{2b-n}(X,M(b))\to H^{2b-n}_{b-n-1,\nr}(X,M(b))).
\]
On the other hand, for any $\alpha\in \ker(H^{2b-n}(X,M(b))\to H^{2b-n}_{b-n-1,\nr}(X,M(b)))$, let $\alpha'$ be the image of $\alpha$ under $H^{2b-1}(X,M(b))\xrightarrow{\rm restr.}H^{2b-n}(F_{b-n}X,M(b))$. Then $\alpha'$ is in the kernel of the restriction map $H^{2b-n}(F_{b-n}X,M(b))\to H^{2b-n}(F_{b-n-1}X,M(b))$, hence by \Cref{lem-a claim} there is an element $\beta \in \CH^b(F_{b-n}X,n;M)$ such that $\cl(\beta)=\alpha'\in H^{2b-n}(F_{b-n}X,M(b))$. Moreover, because the first two columns in (\ref{equ-diagram for main thm2}) are exact and the morphism
$$\varinjlim\limits_{\codim Z=b-n+1}\CH^{n-1}(Z,n-1;M)\xrightarrow{\cl}\varinjlim\limits_{\codim Z=b-n+1}H^{n-1}(Z,M(n-1))$$
is an isomorphism. Therefore, we can lift $\beta$ to be an element of $\CH^b(X,n;M)$, say $\gamma$. Then $(\cl(\gamma)-\alpha)\in H^{2b-n}(X,M(b))$ maps to $\cl(\beta)-\alpha'=0$ in $H^{2b-n}(F_{b-n}X,M(b))$, which means that $\alpha=\cl(\gamma)+\cl(\omega)=\cl(\gamma+\omega)$ for some $\omega\in \CH^b(X,n;M)$ comes from $\varinjlim\limits_{\codim Z=b-n+1}\CH^{n-1}(Z,n;M)\cong \varinjlim\limits_{\codim Z=b-n+1}H^{n-2}(Z,M(n-1))$. Therefore, 
\[
\im(\CH^b(X,n;M)\xrightarrow{\cl} H^{2b-n}(X,M(b)))= \ker(H^{2b-n}(X,M(b))\xrightarrow{\rm restr.} H^{2b-n}_{b-n-1,\nr}(X,M(b))).\qedhere
\]
\end{proof}

Now, we are going to construct a morphism $\theta:H^{2b-n}_{b-n-1,\nr}(X,M(b))\to \CH^b(X,n-1;M)$. 
Let us consider the commutative diagram (\ref{equ-diagram for main thm2}) again, where the second row and the third row are exact by \Cref{lem-a claim} and \Cref{prop-exact at H(X)}. For every $\alpha\in H^{2b-n}_{b-n-1,\nr}(X,M(b))$, we can find an element $\alpha'\in H^{2b-n}(F_{b-n}X,M(b))$ such that $\alpha'$ maps to $\alpha$ under the natural restriction morphism $H^{2b-n}(F_{b-n}X,M(b))\to H^{2b-n}(F_{b-n-1}X,M(b))$.
Note that we have the following commutative diagram 
\begin{equation}\label{equ-construct a connected morphism I}
\xymatrix{H^{2b-n}(F_{b-n}X,M(b))\ar[rrr]^{\partial}\ar[d]_{\iota_*\circ \mathrm{cl}^{-1}\circ \partial}&&&\varinjlim\limits_{\codim Z=b-n+1}H^{n-1}(Z,M(n-1))\\
 \CH^b(X,n-1;M)&&&\varinjlim\limits_{\codim Z=b-n+1}\CH^{n-1}(Z,n-1;M)\ar[u]^{\sim}_{\mathrm{cl}}\ar[lll]_{\iota_*},
       }
\end{equation}
where $Z$ runs through all closed subset of $X$ with $\codim(Z)\coloneqq\dim(X)-\dim(Z)=b-n+1$ and the left vertical arrow is an isomorphism by \Cref{prop-higher cycle class map for any field}. Now, we define 
\begin{equation}\label{equ-construct a connected morphism II}
\theta(\alpha)\coloneqq \iota_*\circ \mathrm{cl}^{-1}\circ \partial(\alpha')\in \CH^b(X,n-1;M), 
\end{equation}
and this is well-defined (i.e., independent of the choice of $\alpha'$) as 
\[
\CH^b(F_{b-n}X,n;M)\to H^{2b-n}(F_{b-n}X,M(b))\to H^{2b-n}(F_{b-n-1}X,M(b))
\]
is exact. Furthermore, this construction is compatible with flat pull-back and proper push-forward as morphisms $\iota_*$, $\cl$ and $\partial$ are.

\begin{thm}\label{thm-main thm2}
Let $X$ be a quasi-projective equi-dimensional scheme defined over a field $k$, and $M\coloneqq \Z/m\Z$, where $m$ is an integer invertible in $k$. For any $b,n\geq 0$, we have the canonical long exact sequence
\begin{small}
\begin{align}\label{equ-main thm2}
    \to\CH^b(X,n;M)\xrightarrow{\mathrm{cl}}H^{2b-n}(X,M(b))\xrightarrow{\rm restr.} H^{2b-n}_{b-n-1,\nr}(X,M(b))\xrightarrow{\theta} \CH^b(X,n-1;M)\to,
\end{align}
\end{small}

where $\CH^b(X,n;M)\coloneqq \CH_{\dim(X)-b}(X,n;M)$ is the Bloch's higher Chow groups with coefficient $M$, and $\mathrm{cl}\colon\CH^b(X,n;M)\to H^{2b-n}(X,M(b))$ is the higher cycle class map defined in \Cref{construction higher cycle class map}. The morphism $\theta:H^{2b-n}_{b-n-1,\nr}(X,M(b))\to \CH^b(X,n;M)$ is constructed in (\ref{equ-construct a connected morphism II}) for $b,n\in \Z$.
\end{thm}
\begin{proof}
By \Cref{cor-GL01} and \Cref{prop-exact at H(X)}, it suffices to show that
\begin{small}
\begin{equation}\label{equ-four terms}
 H^{2b-n}(X,M(b))\xrightarrow{\rm restr.} H^{2b-n}_{b-n-1,\nr}(X,M(b))\xrightarrow{\theta} \CH^b(X,n-1;M) \xrightarrow{\cl} H^{2b-n+1}(X,M(b))
\end{equation}
\end{small}

is exact.

\emph{$\bullet$ (\ref{equ-four terms}) is exact at $H^{2b-n}_{b-n-1,\nr}(X,M(b))$}: This is a complex as the following sequence 
\[
H^{2b-n}(X,M(b))\xrightarrow{\rm restr.} H^{2b-n}(F_{b-n}X,M(b))\xrightarrow{\partial} \varinjlim\limits_{\codim Z=b-n+1}H^{n-1}(Z,M(n-1))
\]
is exact by \Cref{prop-Gysin sequence}. For any $\alpha \in \ker(H^{2b-n}_{b-n-1,\nr}(X,M(b))\xrightarrow{\theta} \CH^b(X,n-1;M))$, we can lift $\alpha$ to an element in $H^{2b-n}(F_{b-n}X,M(b))$, namely $\alpha'$, and $\alpha'$ maps to 
\[
\beta\coloneqq \partial(\alpha') \in \varinjlim\limits_{\codim Z =b-n+1}H^{n-1}(Z,M(n-1))\xlongleftarrow[\sim]{\cl}\varinjlim\limits_{\codim Z=b-n+1}\CH^{n-1}(Z,n-1;M)
\]
by \Cref{prop-higher cycle class map for any field}. Then $\beta$ can be lifted to $\gamma \in \CH^b(F_{b-n}X,n;M)$ with $\partial(\gamma)=\cl^{-1}(\beta)$ as $$0=\theta(\alpha)=\iota_*\circ \cl^{-1}(\beta)\in \CH^b(X,n-1;M).$$
Let $\alpha''\coloneqq \cl(\gamma)$ be the image of $\gamma$ under $\CH^b(F_{b-n}X,n;M)\xrightarrow{\cl} H^{2b-n}(F_{b-n}X,M(b))$, and then $$\partial(\alpha'-\alpha'')=\beta-\partial(\cl (\gamma))=\beta-\cl(\partial(\gamma))=0\in \varinjlim\limits_{\codim Z=b-n+1}H^{n-1}(Z,M(n-1)),$$
which implies $\alpha'-\alpha''$ belongs to the image of the restriction morphism $$H^{2b-n}(X,M(b))\to H^{2b-n}(F_{b-n}X,M(b))$$
by \Cref{prop-Gysin sequence} again. Note that $${\rm restr.}(\alpha'-\alpha'')={\rm restr.}(\alpha')-{\rm restr.}\circ \cl(\gamma)=\alpha\in H^{2b-n}_{b-n-1,\nr}(X,M(b)),$$
so we get $\alpha \in \im(H^{2b-n}(X,M(b))\to H^{2b-n}(F_{b-n}X,M(b)))$.

\emph{$\bullet$ (\ref{equ-four terms}) is exact at $\CH^b(X,n-1;M)$}: This is a complex as the following $$H^{2b-n}(F_{b-n}X,M(b))\xrightarrow{\partial} \varinjlim\limits_{\codim Z=b-n+1}H^{n-1}(Z,M(n-1))\xrightarrow{\iota_*} H^{2b-n+1}(X,M(b))$$ is exact by \Cref{prop-Gysin sequence}. For any $\alpha \in \ker(\CH^b(X,n-1;M)\xrightarrow{\cl} H^{2b-n+1}(X,M(b)))$, $\alpha$ can be lifted to 
\[
\beta \in \varinjlim\limits_{\codim Z=b-n+1}\CH^{n-1}(Z,n-1;M)\cong \varinjlim\limits_{\codim Z=b-n+1}H^{n-1}(Z,M(n-1))
\]
by \Cref{prop-higher cycle class map for any field}. Furthermore, $\cl(\beta)$ can be lifted to $\gamma \in H^{2b-n}(F_{b-n}X,M(b))$ as $$\iota_*\circ\cl(\beta)=\cl(\alpha)=0\in H^{2b-n+1}(X,M(b)),$$ and note that we have 
\begin{equation*}
    \iota_*\circ \cl^{-1}\circ \partial(\gamma)=\iota_*\circ \cl^{-1}\circ \cl(\beta)=\iota_*(\beta)=\alpha\in \CH^b(X,n-1)
\end{equation*}
i.e., $\alpha\in \im(H^{2b-n}_{b-n-1,\nr}(X,M(b))\to \CH^b(X,n-1;M))$ by (\ref{equ-construct a connected morphism I}).
\end{proof}

\begin{cor}\label{cor-kernel of map}
The kernel of the natural map $H^{b+1}(X,M(b))\to H^{b+1}_{\nr}(X,M(b))$ is isomorphic to $\CH^b(X,b-1;M)$.
\end{cor}
We apply the above result to the special case $M=\Z/\ell^r$ for some prime $\ell$ invertible in $k$ and take direct limits over $r$, to obtain the following:
\begin{cor}
    Let $X$ be a quasi-projective equi-dimensional scheme defined over a field $k$, and let $\ell$ be a prime which is not equal to $\mathrm{char}(k)$. For any $b,n\geq 0$, we have the canonical long exact sequence
\begin{small}
 \begin{align*}
\to\CH^b(X,n;\mathbb{Q}_{\ell}/\mathbb{Z}_{\ell})\xrightarrow{\mathrm{cl}}H^{2b-n}(X,\mathbb{Q}_{\ell}/\mathbb{Z}_{\ell}(b))\xrightarrow{\rm restr.} H^{2b-n}_{b-n-1,\nr}(X,\mathbb{Q}_{\ell}/\mathbb{Z}_{\ell}(b))\xrightarrow{\theta} \CH^b(X,n-1;\mathbb{Q}_{\ell}/\mathbb{Z}_{\ell})\to.
\end{align*}   
\end{small}

If $k=\mathbb{C}$, then we can replace $\mathbb{Q}_{\ell}/\mathbb{Z}_{\ell}$ by $\mathbb{Q}/\Z$.
\end{cor}
\begin{rem}
{\rm One thing to remark here is that we only use the `equi-dimensional' and `quasi-projective' assumptions in \Cref{thm-main thm2} for the construction of higher cycle class maps for singular varieties. Even for smooth $X$, we still need the cycle class maps for singular varieties and require it to be at least compatible with localization sequences. Therefore, if one can construct the cycle class map $\CH_i(*,n;\Z/m)\to H_{2i+n}^{BM}(*,\Z/m(i))$ compatible with localization sequences as \Cref{prop-higher cycle class map for any field} for \emph{any} scheme $X$ over the field $k$, where $m$ is invertible in $k$, then (\ref{equ-main thm2}) is also exact for this $X$ and $m$. 
}
\end{rem}
\begin{cor}\label{cor-an exact sequence about torsion cycle}
Let $X$ be a smooth quasi-projective variety defined over an algebraically closed field $k$, and let $\ell$ be a prime which is not equal to $\mathrm{char}(k)$. There is an exact sequence
\begin{equation}\label{equ-an exact sequence about torsion cycle}
    H^{2n-2}_{n-3,\nr}(X,\Q_{\ell}/\Z_{\ell}(n))\to \CH^n(X)[\ell^{\infty}]\xrightarrow{\lambda^n} H^{2n-1}(X,\Q_{\ell}/\Z_{\ell}(n))\to H^{2n-1}_{n-2,\nr}(X,\Q_{\ell}/\Z_{\ell})(n)),
\end{equation}
where $\lambda^n$ is the morphism defined by Bloch in \cite{bloch1979torsion} and the first arrow is the composition of $ H^{2n-2}_{n-3,\nr}(X,\Q_{\ell}/\Z_{\ell}(n))\to \CH^n(X,1;\Q_{\ell}/\Z_{\ell})\to \CH^n(X)[\ell^{\infty}]$. If $k=\C$, we can replace $\mathbb{Q}_{\ell}/\mathbb{Z}_{\ell}$ (resp. $\CH^n(X)[\ell^{\infty}]$) by $\mathbb{Q}/\Z$ (resp. $\CH^n(X)_{\rm tor}$).
\end{cor}
\begin{proof}
    In fact, we have the following commutative diagram
    \begin{equation*}
        \xymatrix{\CH^n(X)[\ell^{\infty}]\ar[r]^{\rho}\ar[dr]^{\lambda^n}&H^{2n}(X,\Z_{\ell}(n))\\
        H^{n-1}_{\rm Zar}(X,\mathcal{H}^n(\Q_{\ell}/\Z_{\ell}(n)))\ar[r]^{\gamma}\ar@{->>}[u]^{\alpha}&H^{2n-1}(X,\Q_{\ell}/\Z_{\ell}(n))\ar[u]^{\beta}\\
        \CH^n(X,1;\Q_{\ell}/\Z_{\ell})\ar[ur]^{\cl}\ar[u]_{\cong}^{\omega}\nospacepunct{,}&
      }
    \end{equation*}
    where
    \begin{itemize}
        \item morphisms $\alpha$, $\beta$, $\gamma$ and $\rho$ are the same as in \cite[Corollary 1]{colliot1983torsion};
        \item the morphism $\lambda^n$ is given by \cite[Section 2]{bloch1979torsion}, and it makes the two triangles above commute by \cite[Theorem 4.3]{colliot1993cycles} and \cite{bloch1979torsion};
        \item the morphism $\omega$ induced by higher cycle class maps is an isomorphism: applying the Gersten's conjecture for higher Chow groups and (\ref{equ-ss for higher chow}), there is a spectral sequence ($m\coloneqq \ell^r$)
        \begin{equation*}
            E_2^{p,q}(n)\coloneqq H^p_{\rm Zar}(X,\mathcal{CH}_X^n(-q))\Rightarrow \CH^n(X,-p-q;\Z/m),
        \end{equation*}
        which is compatible with the spectral sequence
        \begin{equation*}
            E_2^{p,2n+q}\coloneqq H^p_{\rm Zar}(X,\mathcal{H}^{2n+q}_X(n))\Rightarrow H^{2n+p+q}_{\e t}(X,\Z/m(n))
        \end{equation*}
        via higher cycle class maps. We note that $E_2^{p,q}(n)=0$ if $n>-q$ using Gersten's resolution. Moreover, by \Cref{cor-GL01} we have for $n\leq -q$ that $E_2^{p,q}(n)\cong H^p_{\mathrm{Zar}}(X,\mathcal H^{2n+q}_X(\Z/\ell^r\Z(n)))$, which vanishes if $p-q>2n$ again by Gersten's resolution. So the spectral sequence computing $\CH^n(X,1,\Z/\ell^r\Z)$ has only possibly non-zero term $E_2^{n-1,-n}(n)$, thus we get the following   
        \[
        \CH^n(X,1;\Z/\ell^r\Z)\xleftarrow{\sim}E_2^{n-1,-n}(n)=H^{n-1}_{\rm Zar}(X,\mathcal{CH}_X^n(n))\xrightarrow{\sim} H^{n-1}_{\rm Zar}(X,\mathcal{H}^{n}_X(\Z/\ell^r\Z(n))),
        \]
        where the last isomorphism follows from \Cref{cor-GL01} again. Compatibility of the spectral sequences implies that this factors through the higher cycle class map $\CH^n(X,1,\Z/\ell^r\Z)\to H^{2n-1}(X,\Z/\ell^r\Z(n))$ via the edge morphism $$H^{n-1}_{\mathrm{Zar}}(X,\mathcal H^n_X(\Z/\ell^r \Z(n)))\to H^{2n-1}(X,\Z/\ell^r\Z(n)).$$ Taking direct limits of $r$, we get the isomorphism $\omega$.
    \end{itemize}
    Thus, we have the commutative diagram
    \begin{equation*}
        \xymatrix{\CH^n(X)[\ell^{\infty}]\ar[r]^{\lambda^n}&H^{2n-1}(X,\Q_{\ell}/\Z_{\ell}(n))\ar@{=}[d]\\
        \CH^n(X,1;\Q_{\ell}/\Z_{\ell})\ar[r]^{\cl}\ar@{->>}[u]^{\alpha\circ \omega}&H^{2n-1}(X,\Q_{\ell}/\Z_{\ell}(n)),
        }
    \end{equation*}
    and we finish the proof by applying \Cref{thm-main thm2}.
\end{proof}
\begin{rem}
  \rm{A similar result for Griffith groups is proved in \cite[Theorem 1.6(ii)]{schreieder2021refined}, and we prove an analogous exact sequence for non-algebraically fields in \cite{KZ23functorial}. Taking $n=1,\dim(X)$ in (\ref{equ-an exact sequence about torsion cycle}), then one sees that $\lambda^1$ and $\lambda^{\dim(X)}$ are isomorphisms as $H^{2n-2}_{n-3,\nr}(X,\Q_{\ell}/\Z_{\ell}(n))=H^{2n-1}_{n-2,\nr}(X,\Q_{\ell}/\Z_{\ell}(n))=0$ in this case. Taking $n=2$, then we have an exact sequence $$0\to\CH^2(X)[\ell^{\infty}]\xrightarrow{\lambda^2} H^{3}(X,\Q_{\ell}/\Z_{\ell}(2))\to H^{3}_{\nr}(X,\Q_{\ell}/\Z_{\ell}(2)). $$ This covers some results of \cite{bloch1979torsion,colliot1993cycles}. 
  }  
\end{rem}
\section{Another Perspective towards the Long Exact Sequence}
In this section, we conjecture that refined unramified cohomology is a motivic homology theory (i.e., \Cref{conj-truncation for smooth case}). We provide some evidences related to this conjecture, and explain how this is related to known and aforementioned results.
\subsection{A Spectral Sequence Computing Refined Unramified Cohomology}\label{ref unram cohom ss}
In \Cref{properties of homology theory}, we showed that the refined unramified cohomology is a homology theory satisfying the localization sequence (\ref{equ-localization sequence for refined homology}). By \Cref{prop-niveau filtration}, there is thus a spectral sequence 

\begin{equation}\label{equ-niveau for refined}
    ^{\nr}E^{1}_{p,q}\coloneqq \bigoplus_{x\in X_{(p)}}H^{p+q+l,\nr}_{p+q}(x,A(n)) \Rightarrow H^{p+q+l,\nr}_{p+q}(X,A(n))
\end{equation}
for any fixed $l$, where $H_i^{j,\nr}(x,A(n))\coloneqq\varinjlim\limits_{U\subset \overline{\{x\}}}H_i^{j,\nr}(U,A(n))$. Moreover, this spectral sequence is induced by an exact couple 

\begin{equation*}
    \xymatrix{ D\ar[rr]^i&& D\ar[dl]^j\\
&E\ar[ul]^k&
}
\end{equation*}
where $D\coloneqq\bigoplus_{p,q}D^1_{p,q}$, $D^1_{p,q}\coloneqq H^{p+q+l,\nr}_{p+q}(Z_p,A(n))$, and $i\colon D^1_{p,q}\to D^1_{p+1,q-1}$ is the closed embedding. The restriction map $H^*(X,A(n))\to H^*_{*-l,\nr}(X,A(n))$ is induced by the morphism of the associated spectral sequences starting from $^{\nr}E^{p,q}_1$.

We can define $$H^p_{q,\nr}(x,A(n))\coloneqq\varinjlim_{U\subset \overline{\{x\}}}H^p_{q,\nr}(U,A(n))$$ and $$H^p_{q,\nr}(F_jX,A(n))\coloneqq\varinjlim_{F_jX\subset U\subset X} H^p_{q,\nr}(U,A(n)).$$ Then we have the following spectral sequence induced by the coniveau filtration
\begin{equation}\label{equ-coniveau for refined}
    ^{\nr}E_1^{p,q}=\bigoplus_{x\in X^{(p)}}H^{q-p}_{q-l,\nr}(x,A(n-p)) \Rightarrow H^{p+q}_{p+q-l,\nr}(X,A(n)).
\end{equation}
\begin{lem}\label{lem-canonical morphism is an isomorpshim}
There is a canonical morphism $H^p(x,A(n))\to H^p_{q,\nr}(x,A(n))$, which is an isomorphism  for any $q\geq 0$. Moreover, $H^p_{q,\nr}(x,A(n))=0$ if $q<0$.
\end{lem}
\begin{proof}
First, by definition, we have the canonical morphism induced by taking the limit of the restriction morphism $H^p(U,A(n))\to H^p_{q,\nr}(U,A(n))$. Also, we have the canonical morphism
\begin{small}
\[
H^p_{q,\nr}(x,A(n))\coloneqq\varinjlim_{U\subset \overline{\{x\}}}\im(H^p(F_{q+1}U,A(n))\to H^p(F_qU,A(n)))\to \varinjlim_{U\subset \overline{\{x\}}}H^p(F_0U,A(n))\cong H^p(x,A(n)).
\]    
\end{small}

One can see that these both compositions are the identity by definition, so this is an isomorphism. Finally, by definition we have $$H^p_{q,\nr}(x,A(n))=\varinjlim_{U\subset \overline{\{x\}}}H^p_{q,\nr}(U,A(n))=0$$ if $q<0$.
\end{proof}

\begin{rem}\label{rem-remark 5.2}
\rm{ By this lemma, for fixed integer $l\geq 0$, we have $^{\nr}E_1^{p,q}=\bigoplus_{x\in X^{(p)}}H^{q-p}(x,n-p)={E}_1^{p,q}$ if $q\geq l$, otherwise $^{\nr}E_1^{p,q}=0$. Moreover, $^{\nr}E_2^{p,q}={E}_2^{p,q}$ if $q\geq l$, and $^{\nr}E_2^{p,q}=0$ for $q<l$, where $E_1^{p,q}$ is the first page induced by $H_\ast$ from \Cref{prop-niveau filtration}. 
}
\end{rem}
\subsection{Conjecture: Refined Unramified Cohomology is a Motivic Homology Theory}
In this final section of the paper, we propose a conjecture connecting the refined unramified cohomology groups with the hypercohomology of some truncated complex. For this, we let $H^i(-,A(n))$ be the Borel--Moore \'etale homology for $k\neq \C$ and the Borel--Moore singular homology when $k=\C$. It is worth mentioning that during the revision process of this paper, the conjecture was extended and proven by Alexandrou and Schreieder (see \cite{alexandrou2024truncated}). We are grateful for their generosity in sharing their manuscript with us. 

\begin{conj}\label{conj-truncation for smooth case}
Let $A(n)$ be a (twisted) torsion abelian group whose elements have order coprime to the exponential characteristic of $k$, and let $X$ be a smooth variety defined over a field $k$. Consider the identity continuous morphism $\pi\colon X_{\rm \e t}\to X_{\rm Zar}$, and let $\tau^{\geq l}R\pi_*A(n)$ be the $l$-th canonical truncation of the derived complex $R\pi_*A(n)$, where $A(n)$ is the locally constant sheaf with stalk $A(n)$. For any $k,l\geq 0$, there exists a natural isomorphism $f_{k,l}\colon \mathbb{H}^k(X,\tau^{\geq l}R\pi_*A(n))\xrightarrow{\sim} H^k_{k-l,\nr}(X,A(n))$, compatible with localization sequences and making the following diagram commute
\begin{equation*}
    \xymatrix{ \mathbb{H}^k(X,R\pi_*A(n))\ar[r]\ar@{=}[d]&\mathbb{H}^k(X,\tau^{\geq l}R\pi_*A(n))\ar[d]_{\sim}^{f_{k,l}}\\
    H^k(X,A(n))\ar[r]^{\rm{restr.}}&H^k_{k-l,\nr}(X,A(n))
 }
\end{equation*}
Moreover, in case the ground field $k=\mathbb{C}$, we can take $A$ to be any abelian group.
\end{conj}

\begin{rem}
{\rm One can compute directly that 
\[
\mathbb H^k(X,\tau^{\geq l}R\pi_\ast A(n))=\begin{cases}H^k_{\nr}(X,A(n))&\text{ if } l=k;\\ H^k(X,A(n)) &\text{ if } k-l\geq\lceil\frac{k}{2}\rceil,\end{cases}
\]
without \Cref{conj-truncation for smooth case}. Moreover, under the assumption of \Cref{conj-truncation for smooth case}, taking the hypercohomology of the canonical distinguished triangle
\[
R^\ell\pi_\ast A(n)[-l]\to\tau^{\geq l}R\pi_\ast A(n)\to\tau^{\geq l+1}R\pi_\ast A(n)\to R^l\pi_\ast A(n)[-l+1]
\]
would recover \cite[Proposition 7.35]{schreieder2021refined}.
}
\end{rem}

\begin{rem}
\rm{
We also note that assuming \Cref{conj-truncation for smooth case} gives directly an unconditional functoriality along a pullback on refined unramified cohomology. After all if $f\colon X\to Y$ is any morphism between smooth varieties, then we have the following morphism (due to space constraints, we will denote \( A(n) \) simply as \( A \) here.)
\begin{small}
$$\tau^{\geq l}R(\pi_Y)_\ast A\to Rf_*f^*\tau^{\geq l}R(\pi_Y)_\ast A=Rf_*\tau^{\geq l}f^*R(\pi_Y)_\ast A\to Rf_*\tau^{\geq l}R(\pi_X)_*f^*A=Rf_*\tau^{\geq l}R(\pi_X)_*A,$$
\end{small}

where the first arrow is induced by the unit $\id\to Rf_\ast f^\ast$, and the last arrow is induced by the image of unit $\id\to Rf_\ast f^\ast$ under the following composition
\begin{tiny}
 $${\rm Hom}(\id,f_*f^*)\xrightarrow{R(\pi_Y)_*}{\rm Hom}(R(\pi_Y)_*,R(\pi_Y)_*f_*f^*)={\rm Hom}(R(\pi_Y)_*,f_*R(\pi_X)_*f^*)={\rm Hom}(f^*R(\pi_Y)_*,R(\pi_X)_*f^*).$$   
\end{tiny}
Let $p_X$ and $p_X$ be the structure morphism of $X$ and $Y$. Applying $(p_Y)_*$ to the morphism $\tau^{\geq l}R(\pi_Y)_\ast A(n)\to Rf_*\tau^{\geq l}R(\pi_X)A(n)$ and taking $k$-th cohomology, \Cref{conj-truncation for smooth case} then gives the map $f^\ast\colon H^k_{k-l,\nr}(Y,A(n))\to H^k_{k-l,\nr}(X,A(n))$. For the specific construction of the pullbacks of refined unramified cohomology between smooth schemes, we refer to \cite{KZ23functorial}. 
}
\end{rem}
\begin{rem}
\rm{Let us end here by explaining that assuming \Cref{conj-truncation for smooth case} provides another proof of \Cref{thm-main thm2} for smooth $X$. So let $X$ be a smooth variety over a field $k$ and $m$ be an an integer invertible in $k$. As a corollary of Suslin--Voevodsky \cite[Theorem 1.1] {geisser2001bloch} and Voevodsky \cite[Corollary 2]{Voevodsky2002motivic}, there is a natural isomorphism 
\begin{equation}\label{equ-isomorphism for higher chow groups for smooth case}
  \CH^b(X,n;\Z/m)\cong\mathbb H^{2b-n}(X,\tau^{\leq b}R\pi_\ast \mu_m^{\otimes b}).  
\end{equation}
Thus taking the hypercohomology of the canonical distinguished triangle
\[
\tau^{\leq b}R\pi_* \mu_m^{\otimes b}\to R\pi_* \mu_m^{\otimes b}\to \tau^{\geq b+1}R\pi_* \mu_m^{\otimes b}\to \tau^{\leq b}R\pi_* \mu_m^{\otimes b}[1],
\]
and applying \Cref{conj-truncation for smooth case} and (\ref{equ-isomorphism for higher chow groups for smooth case}), we obtain the following long exact sequence again
\begin{align*}
\cdots\to\CH^b(X,n;\Z/m)\xrightarrow{\mathrm{cl}}H_{\e t}^{2b-n}(X, \mu_m^{\otimes b})\to H^{2b-n}_{b-n-1,\nr}(X, \mu_m^{\otimes b})\to \CH^b(X,n-1;\Z/m)\to\cdots.
\end{align*}
}
\end{rem}

\bibliographystyle{alpha}
\bibliography{ref}
\end{document}